%% file: Final_Edits_11_11_21_KO.tex
\documentclass[10pt, a4paper,reqno]{amsart}

\usepackage[total={6.75in, 8in}]{geometry}
\usepackage{graphicx}
\usepackage{booktabs}
\usepackage{float}
\usepackage{esint}
\usepackage{cancel}
\usepackage[usenames,dvipsnames]{xcolor}

\usepackage[]{hyperref}
\hypersetup{
    backref         =   true,       
    pagebackref     =   true,               
    hyperindex      =   true,                
    colorlinks      =   true,                
    breaklinks      =   true,                
    urlcolor        =   black,                
    linkcolor       =   WildStrawberry,                
    bookmarks       =   true,                 
    bookmarksopen   =   false,
    filecolor       =   black,
    citecolor       =   WildStrawberry,
    linkbordercolor =   blue
}

\usepackage{cleveref}

\title{
    A Novel Method for Generating 3D Hydrodynamic Conservation Laws}
\author{James Wing Chee Graham}
\address{Western University, London, Ontario}
\author{Katie Oliveras}
\address{Seattle University, Seattle, Washington}
\author{Olga Trichtchenko}
\address{Western University, London, Ontario}

\date{November 28, 2021}

\begin{document}

\begin{abstract}
    We extend the recent works of Oliveras \cite{oliveras2020weak} to develop a new nonlocal formulation of the water-wave problem for a three-dimensional fluid with a two-dimensional free surface for an inviscid and irrotational fluid over a flat bathymetry. Using this new formulation, we show how one can systematically derive Benjamin \& Olver's twelve conservation laws without explicitly relying on the underlying Lie symmetries. This allows us to make draw new conclusions about conservation laws and posit the potential existence of additional, nonlocal, conservation laws for the water-wave problem.
\end{abstract}
\maketitle

\section{Introduction}
    In 1982, Benjamin and Olver derived the Lie symmetry groups for an inviscid, irrotational fluid in two- and three-dimensions with a free-boundary one- and two-dimensional free-boundary \(\eta\) \cite{benjamin1982hamiltonian}. From these symmetries, the conserved densities were derived in terms of quantities evaluated on the boundary of the fluid domain. Olver subsequently argued that for the 2D problem, the eight conservation laws that arise from the 2D Lie symmetry group are unique and exhaustive as far as deriving conservation laws via infinitesimal symmetries for the free boundary problem. However, a rigorous proof of their exhaustiveness has yet to be seen. 

    To the best of our knowledge, Benjamin and Olver's work represents the only systematic derivation of the appropriate conserved densities for the water wave problem using the underlying Lie symmetries. In \cite{oliveras2021different}, Oliveras and Calatola-Young presented a new systematic derivation of the conserved densities of the 2D water wave problem via a weak formulation, which we seek to extend to fluids in three spatial dimenions.  

    We consider the nonlocal/nonlocal formulation presented in \cite{oliveras2020weak} for an inviscid, irrotational, incompressible fluid posed on the whole line with a two-dimensional free surface. We then illustrate a process for deriving Benjamin and Olver’s conservation laws without explicitly using the symmetries. 

    The water-wave problem (or hydrodynamic problem) in three spatial dimensions is known to possess twelve conservation laws as found by Benjamin \& Olver \cite{benjamin1982hamiltonian}. We provide here an elucidation of the techniques for deriving the twelve conserved densities of the 3D problem here, generalizing the methods of \cite{oliveras2020weak} to the case of waves in three dimensions with a two-dimensional free surface. 

    In \Cref{sec:prelim}, we present the equations of motion for an inviscid, irrotational fluid in three-dimensions with a free surface and define conserved density in terms of boundary variables.  In \Cref{sec:reform}, we introduce a generalization of the non-local integral formulation that will be the foundation for deriving conservation laws.  In \Cref{sec:densities}, we provide full and extensive derivations of the conserved densities of this problem. Finally, in \Cref{sec:surfTen}, we include the effects of surface tension, and briefly discuss the implications for our conserved density results. 

\section{Preliminaries}\label{sec:prelim}
\begin{figure}[hbt!] 
        \centering
        \includegraphics[width=.55\textwidth]{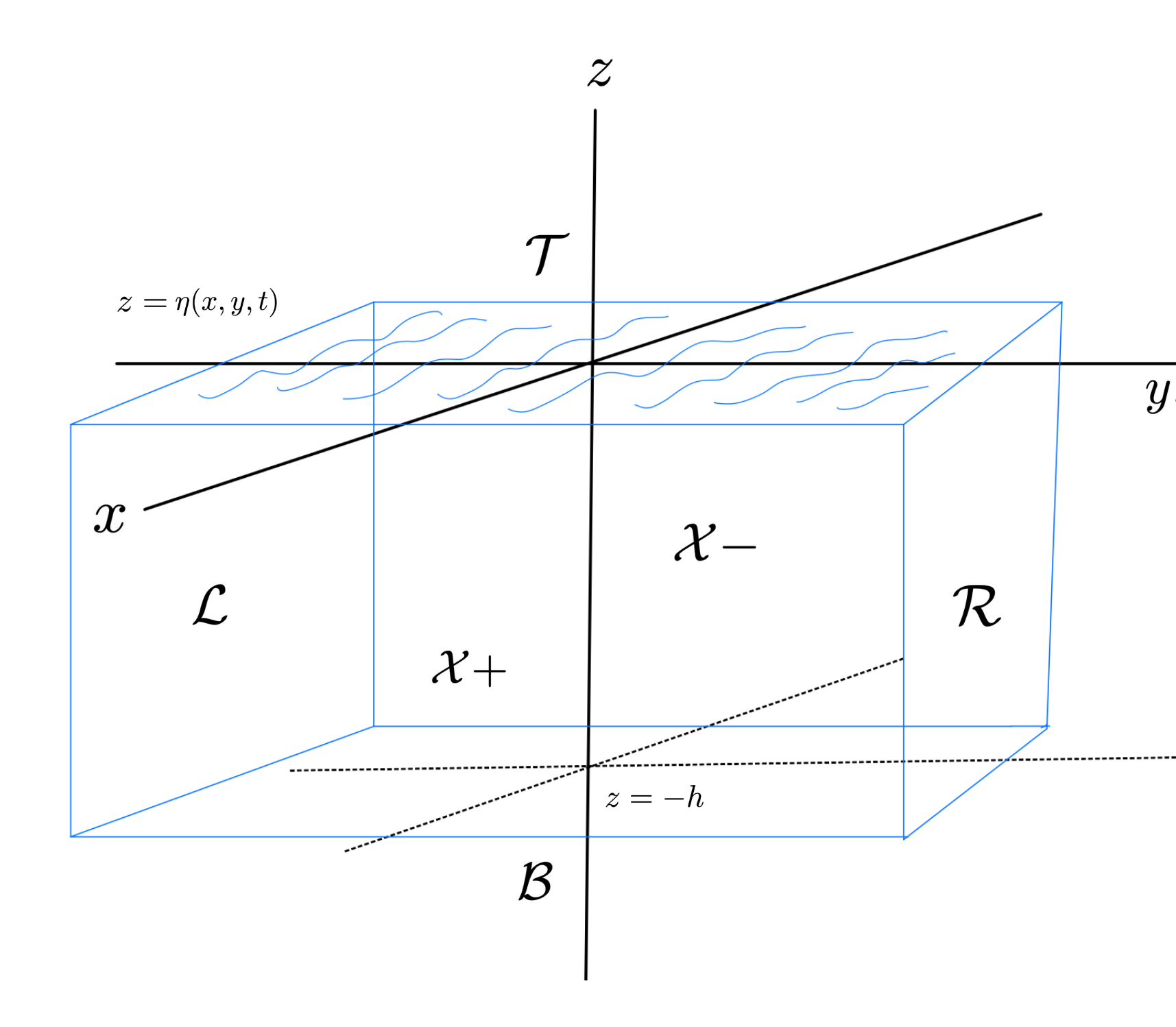}
        \caption{The fluid domain for the three-dimensional water-wave problem.}\label{fig:fluidDomain}
    \end{figure}

    We being by considering the equations of motion for a three-dimensional, irrotational, and inviscid fluid in finite depth that can be described as follows:
    \begin{subequations}\label{hydro_3d}
        \begin{align} \label{hydro1_3d}
            \phi_{xx} + \phi_{yy} + \phi_{zz} = 0,  &&  (x,y,z) \in \mathcal{V},\\ \label{hydro2_3d}
            \phi_{t} + \frac{1}{2}|\nabla\phi|^{2} + gz + p/\rho = 0,  && (x,y,z) \in \mathcal{V},  \\ \label{hydro3_3d}
            \phi_{z} = 0,   && z = -h, \\ \label{hydro4_3d}
            \eta_{t} + \phi_{x}\eta_{x} + \phi_{y}\eta_{y} = \phi_{z},  && z = \eta(x,y,t), \\ \label{hydro5_3d}
            \phi_{t} + \frac{1}{2}|\nabla\phi|^{2} + g\eta = 0,   && z = \eta(x,y,t) ,
        \end{align}
    \end{subequations} 
    where \(\phi(x,y,z,t)\) is the velocity potential in the bulk of the fluid and \(\eta(x,y,t)\) is the two-dimensional free-surface surface evolving in time. For simplicity, we restrict ourselves to a free-surface \(z = \eta(x,y,t)\) that is single-valued, continuously differentiable, and decays sufficiently fast as \(|x|,|y|\infty\). We also impose that our velocity potential \(\phi\to 0\) sufficiently fast as \(|x|,|y|\to\infty\) thus enforcing that \(|\nabla \phi|\to 0\) sufficiently fast as \(|x|,|y| \to \infty\). See \Cref{fig:fluidDomain} for reference.

    The 3D problem involves flux integrals over six surfaces (\(\mathcal{T}, \mathcal{B}, X+, X-, L, R\)) defined in \Cref{fig:fluidDomain} that encloses a fluid volume \(V\). We will denote the boundary of this volume \(V\) as \(\partial \mathcal{V}\), and use \(\partial \mathcal{V}\) to denote the union of all 6 surfaces.

    \subsection{Conserved Densities for the Water-Wave Problem}
        As an objective of this paper is to define conservation laws for the water-wave problem, we begin by stating the definition.  
 
        Following the work of \cite{benjamin1982hamiltonian,olver1983conservation}, we define a scalar function function T depending on x, y, t, \(\eta\), \(\phi\), and the derivatives of \(\eta\) and \(\phi\) over the free-surface \(\eta\) to be a \emph{conserved density} for the free-boundary problem if there exists a vector function
        \begin{equation} 
            \vec{F} = \begin{pmatrix} F_{1} \\ F_{2} \\ F_{3} \end{pmatrix},
        \end{equation} 
        with \(F_{1}, F_{2}, F_{3}\) depending on \(x, y, z, t\), and \(\phi\) (as well as the derivatives of \(\phi\)) in the region \(V\), and a function \(W\) depending on \(x, y, t, z, \phi\) (as well as their derivatives) evaluated along the free-surface such that for all solutions to the of the free boundary problem,

        \begin{displaymath}
            D_{t}T = \vec{F} \cdot \vec{n} + D_{x}W, \text{ on } \mathcal{T}, \text{ with }\qquad \nabla \cdot \vec{F} = 0, \quad \text{in } V.
        \end{displaymath}

        In \cite{benjamin1982hamiltonian}, Benjamin \& Olver find the twelve conserved densities of the 3D hydrodynamic problem shown in \Cref{densities}.

        \begin{table}[H]
            \centering
            \caption{Table of Conserved Densities from Benjamin \& Olver \cite{benjamin1982hamiltonian}.}\label{densities}
            \renewcommand*\arraystretch{2}
            \begin{tabular}{lcl}
                \toprule
                \(\displaystyle \boldsymbol{T_{1}} = -q\eta_{x}\) && (horizontal   impulse) \\
                \(\displaystyle \boldsymbol{T_{2}} = -q\eta_{y}\) && (horizontal    impulse) \\
                \(\displaystyle \boldsymbol{T_{3}} = \mathcal{H} = \frac{1}{2}q\eta_{t} + \frac{1}{2}g\eta^{2} \)&& (energy) \\
                \(\displaystyle \boldsymbol{T_{4}} = \eta \)&& (mass) \\
                \(\displaystyle \boldsymbol{T_{5}} = q + gt\eta\) && (vertical    impulse) \\
                \(\displaystyle \boldsymbol{T_{6}} = x\eta + tq\eta_{x}\) && (horizontal    \(x\)   Gallilean   boost) \\
                \(\displaystyle \boldsymbol{T_{7}} = y\eta + tq\eta_{y}\) && (horizontal    \(y\)   Gallilean   boost) \\
                \(\displaystyle \boldsymbol{T_{8}} = \frac{1}{2}\eta^{2} - tT_{5} + \frac{1}{2}gt^{2}\eta\) && (vertical   Gallilean   boost) \\
                \(\displaystyle \boldsymbol{T_{9}} = (x\eta_{y} - y\eta_{x})q\) && (horizontal   rotation) \\
                \(\displaystyle \boldsymbol{T_{10}} = (x + \eta\eta_{x})q + gt(x\eta + tq\eta_{x}) - g \frac{1}{2} t^{2} q \eta_{x}\) && (\(x\)   gravity   compensated   rotation) \\ 
                \(\displaystyle \boldsymbol{T_{11}} = (y + \eta\eta_{y})q + gt(y\eta + tq\eta_{y}) - g \frac{1}{2} t^{2} q \eta_{y}\) && (\(y\)   gravity   compensated   rotation) \\
                \(\displaystyle \boldsymbol{T_{12}} = (\eta - x\eta_{x} - y\eta_{y})q + t(9gT_{8} - 5 T_{3}) + \frac{9}{2}gt^{2}T_{5} - \frac{3}{2}g^{2}t^{3}T_{4}\) && (vertical   acceleration) \\
                \bottomrule
            \end{tabular}
        \end{table}

        With this definition, a generalized conservation law in surface variables takes the following form:
        \begin{equation} \label{genconsvOriginal}
            \frac{d}{dt} \iint_{\mathcal{T}} \; [T] \;  \,dx\,dy = - \iint_{\Gamma} \vec{F} \cdot \vec{n} \; dS + W \: \large{|}_{\partial \mathcal{T}},
        \end{equation} 
        
        where \(T\) represents the generalized conserved density, \(\vec{F}\) is a generalized force vector field, \(\vec{n}\) is the outward unit normal to the surface element \(dS\), \(\Gamma = \partial \mathcal{V} \backslash \mathcal{T}\), and \(W\) is a generalized work flux over the free surface \(\mathcal{T}\). 
        In \eqref{genconsvOriginal}, we have introduced the following notation that will be used throughout the paper: 
        \begin{equation}\label{doubleiterated1}
            \iint_{\mathcal{T}} f(x,y,z,t) \,dx\,dy = \iint_{\mathcal{T}} f(x,y,\eta(x,y,t),t) \,dx\,dy = \int_{-\infty}^{\infty}\int_{-\infty}^{\infty}  f(x,y,\eta(x,y,t),t) \,dx\,dy.
        \end{equation}
        Likewise, we define the integral along the bottom boundary to be 
        \begin{equation}
            \iint_{\mathcal{B}} f(x,y,z,t) \,dx\,dy = \iint_{\mathcal{B}} f(x,y,-h,t) \,dx\,dy = \int_{-\infty}^{\infty}\int_{-\infty}^{\infty} f(x,y,-h,t) \,dx\,dy.
        \end{equation} 
        Equation (\ref{genconsvOriginal}) expresses the fact that the rate of changes in total density equals the sum of the fluxes over the fixed boundary \(\Gamma\) and the boundary of the free surface. It can be shown that an equivalent form of such a generalized conservation law can be written as
        \begin{equation}\label{genconsvOriginal2}
            \frac{d}{dt} \iint_{\mathcal{T}} \; [T] \;  \,dx\,dy = - \iint_{\mathcal{B}} U  \,dx\,dy +  \iint_{\mathcal{T}} [ \frac{d}{dx} V]  \,dx\,dy +  \iint_{\mathcal{T}} [ \frac{d}{dy} V]  \,dx\,dy.
        \end{equation}

\section{A Reformulation of the equations of Motion}\label{sec:reform}

    \subsection{The formulation of Ablowitz, Fokas, and Musslimani}

        The Ablowitz', Fokas, and Musslimani (AFM) formulation of water waves is one of several surface-variable formulations of the water-wave problem \cite{ablowitz2006new}.  By letting \(q(x,y,t) = \phi(x,y,\eta(x,y,t),t)\) a \(\eta(x,y,t) = z\vert_{\mathcal{T}}\), the authors formulate \eqref{hydro_3d} as a system of two equations for the unknowns surface variables \(q(x,y,t)\) and \(\eta(x,y,t)\) as 
        \begin{subequations}\label{eqns:AFM}
            \begin{eqnarray}
                &&\label{local}\displaystyle q_{t} + \frac{1}{2}|\nabla^{\perp} q|^{2} + g\eta - \frac{1}{2}\frac{(\eta_{t}+ \nabla^{\perp} q \cdot \nabla^{\perp} \eta)^{2}}{1 + |\nabla \eta|^{2}} = 0, \\
                &&\label{nonlocal}\iint_{\mathcal{T}}e^{i(k_{x}x + k_{y}y)} \left[i\eta_{t}\cosh{(\vert\vec k\vert(z + h))} + \frac{(\vec{k} \cdot \nabla^{\perp} q )}{\vert\vec k\vert}\sinh{(\vert\vec k\vert(z + h))}\right] \,dx\,dy = 0,\qquad \forall ~\vec{k}\in\mathbb{R}^2,
            \end{eqnarray}
        \end{subequations} 
        where \(\nabla^{\perp} = (\partial_x, \partial_y)\).  In what follows, we will refer to \eqref{eqns:AFM} as the Local/Nonlocal formulation. 

    \subsection{A Nonlocal/Nonlocal Formulation}
        In order to systematically derive the conservation laws from the water-wave problem, we find it useful to reformulate the problem as a system of nonlocal equations in terms of the same surface variables used by \cite{ablowitz2006new}.  Following \cite{ablowitz2006new} and \cite{oliveras2020weak}, we consider the velocity potential \(\phi\) that solves \eqref{hydro_3d} in \(V\) as well as a harmonic test-function \(\psi\) defined in \(V\).  Since \(\psi\) is harmonic, \(\psi_z\) is also harmonic.  Thus, with both \(\phi\) and \(\psi_z\) harmonic in \(V\), via Green's second identity, we have the following integral forms:
        \begin{subequations}
        \begin{align} \label{green1_3d}
        \oiint_{\partial \mathcal{V}} \left[ \psi_{z}(\nabla \phi \cdot \vec{n}) - \phi(\nabla \psi_{z} \cdot \vec{n}) \right] dS = 0,\\ \label{green2_3d}
        \frac{d}{dt} \oiint_{\partial \mathcal{V}} \left[ \psi_{z}(\nabla \phi \cdot \vec{n}) -  \phi(\nabla \psi_{z} \cdot \vec{n}) \right] dS = 0. 
        \end{align}
        \end{subequations}
        These integral forms will lead ultimately to our integral equations of motion determining the twelve conservation laws for this system. It should be noted that upon choosing the harmonic function \(\psi = e^{-i(k_x x + k_y y)}\sinh(\vert k\vert(z + h))\), \eqref{green1_3d} becomes the non-local equation of motion in the AFM formulation stated above in (\ref{nonlocal}) provided that both \(\eta(x,y,t)\) and \(\vert\nabla\phi\vert\) decay to zero sufficiently fast as \(\vert x\vert, \,\vert y \vert \to \infty\) such that the integral makes sense.  

        We can rewrite \eqref{green1_3d} by substituting the boundary conditions given in \eqref{hydro5_3d} and by emposing suitable decay such that the following expression is valid
        \begin{displaymath}
            \iint_{\mathcal{T}} \frac{d}{dt} \left[ \psi \right]  \,dx\,dy - \iint_{\mathcal{T}} \left[  q_{x} \psi_{x}  +   q_{y}\psi_{y} \right] \,dx\,dy =   \iint_{\Gamma} \left[  \phi(\nabla \psi_{z} \cdot \vec{n}) \right] dS.
        \end{displaymath}
        where \(\Gamma = \partial \mathcal{V} \backslash \mathcal{T}\).  Here, we have expressed our equation as a time derivative of an iterated integral evaluated at the free surface over the \(x,y\)-plane equated to flux terms out of the bounding surfaces of our volume \(V\), plus an extra term over the free-surface which will be instrumental in deriving conservation laws. If we use decay arguments to reduce the \(\Gamma\) term to the bottom domain, we find that it actually vanishes due to the bed condition, and that our final integral equation of motion is: 
        \begin{equation}
            \iint_{\mathcal{T}}\frac{d}{dt}  \left[ \psi \right]  \,dx\,dy = \iint_{\mathcal{T}} \left[  q_{x} \psi_{x}  +  q_{y}\psi_{y} \right] \,dx\,dy  - \iint_{\mathcal{B}} \left[  \phi\psi_{zz} \right] \,dx\,dy. \label{eqn:firstInt}
        \end{equation}
        Thus, \Cref{eqn:firstInt} represents our first integral equation of motion for deriving our conservation laws and can be considered a generalization of \eqref{nonlocal}.

        To find our second nonlocal equation of motion, we repeat the same approach using \eqref{green2_3d} as a starting point.  Using Reynold's transport theorem from continuum mechanics,  \eqref{green2_3d} can be written as 
        \begin{align*}
            \frac{d}{dt}\iint_{\mathcal{T}} \left[ \psi_{z} \eta_{t} \right] \,dx\,dy =  \iiint_{ \mathcal{V}}\left[   \nabla \phi_{t} \cdot \nabla \psi_{z}  \right] dV + \oiint_{ \partial \mathcal{V}} \left[ ( \nabla \phi \cdot \nabla \psi_{z}) \nabla \phi \cdot \vec{n} \right] dS.
        \end{align*} 

        Substituting the boundary conditions from \eqref{hydro5_3d} (and after some tedious calculations), we arrive at the expression 
        \begin{align*}
        \frac{d}{dt}\iint_{\mathcal{T}} \left[ \psi_{z} \eta_{t} \right] \,dx\,dy =  \iint_{ \mathcal{T}}\left[  q_{t}  \nabla \psi_{z} +  ( \nabla^{\perp} q \cdot \nabla^{\perp} \psi_{z} ) \nabla \phi \right] \cdot \vec{n} \,dS + \iint_{\Gamma}\left[\phi_{t} \nabla \psi_{z}  \right]  \cdot \vec{n} \,dS,
        \end{align*}
        where \(\Gamma = \partial \mathcal{V} \backslash \mathcal{T}\). As before, if we impose the appropriate decay conditions, our second integral equation of motion is given by 
        \begin{equation}
        \frac{d}{dt}\iint_{\mathcal{T}} \left[ \psi_{z} \eta_{t} \right] \,dx\,dy -  \iint_{ \mathcal{T}}\left[  q_{t}  \nabla \psi_{z} +  ( \nabla^{\perp} q \cdot \nabla^{\perp} \psi_{z} ) \nabla \phi \right] \cdot \vec{n} \,dS=  - \iint_{\mathcal{B}}\left[\phi_{t} \psi_{zz} \right] \,dx\,dy. \label{eqn:secondInt}
        \end{equation}

        Thus, Equations \eqref{eqn:firstInt} and \eqref{eqn:secondInt} with \(q_t\) replaced using \eqref{local}, represent our Nonlocal/Nonlocal formulation of the water-wave problem.

\section{Conserved Densities} \label{sec:densities}
    Using the nonlocal/nonlocal formulation given by \eqref{eqn:firstInt} and \eqref{eqn:secondInt}, we aim to systematically derive the conserved densities found in \cite{benjamin1982hamiltonian} by varying the harmonic test-function \(\psi\).  Before we begin this task, we first review the Benjamin \& Olver's definition for \emph{conserved density} and briefly summarize the results for a 3d fluid.  

    We define three separate 3D harmonic test functions \(\psi_{n}^{x}\), \(\psi_{n}^{y}\), and \(\psi_{n}^{xy}\) such that:

    \begin{equation} 
        \label{testfunction3d_2} 
        \psi_{n}^{x}(x,z) = \frac{(x + iz)^{n}}{n!}, \qquad \psi_{n}^{y}(y,z) = \frac{(y + iz)^{n}}{n!}, \qquad \psi_{n}^{xy}(x,y) =  \frac{(x + iy)^{n}}{n!}.
    \end{equation}

    We will now present derivations of the twelve conservation laws for the 3D Hydrodynamic problem. The basic idea is to consider various harmonic test functions \(\psi\) such as those defined above.  Substituting these choices for \(\psi\) directly into \eqref{eqn:firstInt} and \eqref{eqn:secondInt} yield all twelve conservation laws. The intuition for this method lies in the fact that conservation of mass falls directly out of the non-local equation in the AFM formulation \cite{ablowitz2006new}. 
        
    Rather than proceeding in the order of conserved densities as presented in \cite{benjamin1982hamiltonian}, we proceed by increasing the value of \(n\) in our harmonic test functions \(\psi_n\).  For the sake of clarity, we present the results here in brief.  The detailed compuations can be found in \Cref{app:A}.
    
    To find our first conservation law, we consider the simplest form of \(\psi\) given by \(\psi_{1}^{x} = x + iz\).  Substituting this choice into (\ref{eqn:firstInt}), yields 
    \begin{displaymath}
        \iint_{\mathcal{T}}\frac{d}{dt}  \left[ (x + iz) \right]  \,dx\,dy - \iint_{\mathcal{T}} \left[  q_{x} \right] \,dx\,dy =  + \iint_{\Gamma} \left[  \phi(\nabla (i) \cdot \vec{n}) \right] dS
    \end{displaymath}
    Simplifying the above expression yields the expression
    \begin{displaymath}
        \frac{d}{dt} \iint_{\mathcal{T}} \left[ \eta \right]  = 0,
    \end{displaymath}
    precisely \(T_4\) as seen in \Cref{densities}. We can proceed in a systematic manner to generate additional conservation laws by increasing the value of \(n\), and separating real and imaginary parts.  The full details are provided in \Cref{app:A}.  

    The only exception to this process is the derivation of \(T_{12}\).  In order to determine this conservation law, we have to consider a test function which is a linear combination \(\psi_3^x\) and \(\psi_3^y\).  To be precise, we take \(\psi = \psi_3^x + \psi_3^y\).  

\section{Surface Tension}\label{sec:surfTen}
        
        Upon including surface tension in the formulation, the hydrodynamic problem in 3D becomes: \\
        
        \begin{subequations}
        \begin{align} 
        &\phi_{xx} + \phi_{yy} + \phi_{zz} = 0  &&  (x,y,z) \in \mathcal{V} \\ 
        &\phi_{t} + \frac{1}{2}|\nabla\phi|^{2} + gz + p/\rho = 0  && (x,y,z) \in \mathcal{V}   \\ 
        &\phi_{z} = 0   && z = -h \\ 
        &\eta_{t} + \phi_{x}\eta_{x} + \phi_{y}\eta_{y} = \phi_{z}  && z = \eta(x,y,t) \\ 
        &\phi_{t} + \frac{1}{2}|\nabla\phi|^{2} + g\eta = \frac{\sigma}{\rho}\nabla \cdot \left[ \frac{\nabla  \eta}{\sqrt{1 + |\nabla \eta|^{2}}} \right]  && z = \eta(x,y,t) 
        \end{align}
        \end{subequations} 
        where \(\sigma\) represents the coefficient of surface-tension. The resulting non-local formulation still holds:        
        \begin{align}
        \frac{d}{dt} \iint_{\mathcal{T}} \left[ \psi \right]  \,dx\,dy - \iint_{\mathcal{T}} \left[  q_{x} \psi_{x}  +   q_{y}\psi_{y} \right] \,dx\,dy =  + \iint_{\Gamma} \left[  \phi(\nabla \psi_{z} \cdot \vec{n}) \right] dS  \\ 
         \frac{d}{dt}\iint_{\mathcal{T}} \left[ \psi_{z} \eta_{t} \right] \,dx\,dy - \iint_{ \mathcal{T}}\left[  q_{t}  \nabla \psi_{z} +  ( \nabla^{\perp} q \cdot \nabla^{\perp} \psi_{z} ) \nabla \phi \right] \cdot \vec{n}\, dS=   \iint_{\Gamma}\left[\phi_{t} \nabla \psi_{z}   \right]  \cdot \vec{n}\, dS 
        \end{align} 
        where in the second equation, there is an implicit dependence on the surface-tension-modified local AFM equation \cite{ablowitz2006new} given by        
        \begin{align}
        q_{t} + \frac{1}{2}|\nabla q|^{2} + g\eta + \frac{(\eta_{t} + \nabla q \cdot \nabla \eta )^{2}}{2(1+|\nabla \eta|^{2})}= \frac{\sigma}{\rho}\nabla \cdot \left[ \frac{\nabla  \eta}{\sqrt{1 + |\nabla \eta|^{2}}} \right]
        \end{align} 
        
        The results for conserved densities of the 3D hydrodynamic problem are almost exactly the same as \Cref{densities} with the inclusion of surface tension, with the exceptions of \(T_{3}\) and \(T_{12}\). The Hamiltonian density \(T_{3}\) contains an an additional term as shown below        
        \begin{align}
        T_{3} = \frac{1}{2}g\eta^{2} + \frac{1}{2}q\eta_{t} + \sigma(\sqrt{1 + |\nabla \eta|^{2}} - 1).
        \end{align} 
        However, with surface tension included, \(T_{12}\) actually ceases to be a conserved density. When attempting to recover \(T_{12}\) in the derivation, one encounters the term:         
        \begin{align}
        \frac{d}{dt} \iint_{\mathcal{T}} \left[ \sigma \sqrt{ 1 + |\nabla \eta|^{2}} \right] \,dx\,dy \neq 0.
        \end{align} 
        
        The fact that \(T_{12}\) is no longer conserved when surface tension is present makes physical sense considering that there is no reason that the arc-length of the free-surface should be time invariant as the system evolves.
        
        \section{Conclusion}
        
        In summary, we have presented a generalization of the non-local AFM formulation for the water wave problem in 3D. We have derived two integral equations of motion that condense Euler's equations for inviscid, irrotational, incompressible water waves into two non-local equations that ultimately yield Olver's twelve conserved densities.     

        We wrote down 3D generalizations of the harmonic test functions, and used them in the twelve derivations of the conservation laws. We omitted laws 2 and 11 as they provided no further useful information due to their symmetry to laws 1 and 10, respectively.        
        
        Modifications to the formulation were also presented for the inclusion of a surface tension term for the interface, and changes to the conserved density results were discussed, noting that \(T_{12}\) no longer remained conserved, and the Hamiltonian density \(T_{3}\) contained an extra surface tension term.         
        
        Further work may include writing this formulation in the travelling waves frame, condensing the two integral equations of motion into one, eliminating the velocity potential variable and ultimately yielding a single equation of motion that can be numerically solved for the free-surface, similar to the work of Oliveras \& Vasan \cite{oliveras2013new}, and Oliveras \cite{oliveras2020weak}.

\bibliographystyle{plain}
\bibliography{references.bib}

\input{appendix.tex}

\end{document}

%% file: appendix.tex
\appendix
\section{Derivation} \label{app:A}
    This appendix contains the detailed computations that give rise to the conservation laws found in \Cref{sec:densities}.
        
    \subsection{Conservation Law 1: \(-q\eta_{x}\)}
        
            We use integral equation of motion (\ref{eqn:secondInt}) and isolate the imaginary part of \(\psi_{2}^{x}\) here to obtain the conservation law.  
            
            \begin{displaymath}
             \frac{d}{dt}\iint_{\mathcal{T}} \Big[ (-\eta + ix) \eta_{t} \Big] \,dx\,dy - \iint_{ \mathcal{T}}\Big[  q_{t}  \nabla (-\eta + ix)+  ( \nabla^{\perp} q \cdot \nabla^{\perp} (-\eta + ix) ) \nabla \phi \Big] \cdot \vec{n}\, dS=  \iint_{\mathcal{B}}\Big[\phi_{t}  \Big]  \,dx\,dy 
            \end{displaymath}
            
            Isolating imaginary part:
            
            \begin{displaymath}
             \frac{d}{dt}\iint_{\mathcal{T}} \Big[ ix \eta_{t} \Big] \,dx\,dy - \iint_{ \mathcal{T}}\Big[  q_{t}  \nabla (ix)+  ( \nabla^{\perp} q \cdot \nabla^{\perp} (ix) ) \nabla \phi \Big] \cdot \vec{n}\, dS = 0
            \end{displaymath}
            
            Expanding and dropping \(i\)'s  
            
            \begin{displaymath}
             \frac{d}{dt}\iint_{\mathcal{T}} \Big[ x \eta_{t} \Big] \,dx\,dy - \iint_{ \mathcal{T}}\Big[  q_{t}(-\eta_{x}) +  q_{x}\eta_{t} \Big]  \,dx\,dy = 0
            \end{displaymath}

            \begin{displaymath}
             \frac{d}{dt}\iint_{\mathcal{T}} \Big[ x \eta_{t} \Big] \,dx\,dy + \iint_{ \mathcal{T}}\Big[  q_{t}\eta_{x} -  q_{x}\eta_{t} \Big]  \,dx\,dy = 0
            \end{displaymath}
            
            Expanding \(q_{t} = \phi_{t} + \phi_{z}\eta_{t}\) and \(q_{x} = \phi_{x} + \phi_{z}\eta_{x}\), factoring out a -1:
            
            \begin{displaymath}
             \frac{d}{dt}\iint_{\mathcal{T}} \Big[ x \eta_{t} \Big] \,dx\,dy - \iint_{ \mathcal{T}}\Big[    \phi_{x}\eta_{t} - \phi_{t}\eta_{x}  \Big]  \,dx\,dy = 0
            \end{displaymath}
            
            Pushing to the boundary, noting the \(\Gamma\) contribution vanishes by decay arguments and constant bathymetry:
            
            \begin{displaymath}
             \frac{d}{dt}\oiint_{\partial \mathcal{V}} \Big[ x \nabla \phi  \Big] \cdot \vec{n}\, dS  -\oiint_{\partial \mathcal{V}} \Big[ \phi_{x} \nabla \phi \Big] \cdot \vec{n}\, dS  - \oiint_{ \partial \mathcal{V}}\Big[ (\phi_{t},0,0) \Big] \cdot \vec{n}\, dS + \cancelto{0}{\iint_{ \Gamma }\Big[ (\phi_{t},0,0) \Big] \cdot \vec{n} \,dS } = 0
            \end{displaymath}
            
            By Green's second identity and the harmonicity of both x and \(\phi\), we can re-write the left most term as
            
            \begin{displaymath}
             \frac{d}{dt}\oiint_{\partial \mathcal{V}} \Big[ \phi \nabla x \Big] \cdot \vec{n}\, dS  -\oiint_{\partial \mathcal{V}} \Big[ \phi_{x} \nabla \phi \Big] \cdot \vec{n}\, dS  - \oiint_{ \partial \mathcal{V}}\Big[ (\phi_{t},0,0) \Big] \cdot \vec{n}\, dS = 0
            \end{displaymath}
            
            Pushing to the surface, we obtain
            
            \begin{displaymath}
             \frac{d}{dt}\iint_{\mathcal{T}} \Big[ -q\eta_{x} \Big] \,dx\,dy  -\oiint_{\partial \mathcal{V}} \Big[ \phi_{x} \nabla \phi \Big] \cdot \vec{n}\, dS  - \oiint_{ \partial \mathcal{V}}\Big[ (\phi_{t},0,0) \Big] \cdot \vec{n}\, dS = 0
            \end{displaymath}
            
            Pushing the two other terms to the bulk, using Green's first identity and the Divergence Theorem, for the first and second terms, respectively.
            
            \begin{displaymath}
             \frac{d}{dt}\iint_{\mathcal{T}} \Big[ -q\eta_{x} \Big] \,dx\,dy  -\oiint_{\mathcal{V}} \Big[ \nabla \phi_{x} \cdot \nabla \phi  + \nabla \cdot (\phi_{t},0,0) \Big]  dV = 0
            \end{displaymath}
            
            Note that \(\nabla \phi_{x} \cdot \nabla \phi =  \frac{1}{2}|\nabla \phi|^{2}_{x} \):
            
            \begin{displaymath}
             \frac{d}{dt}\iint_{\mathcal{T}} \Big[ -q\eta_{x} \Big] \,dx\,dy  -\oiint_{\mathcal{V}} \Big[ \frac{1}{2}|\nabla \phi|^{2}_{x}  + \phi_{tx}  \Big]  dV = 0
            \end{displaymath}
            
            Using dynamic boundary condition:

            \begin{displaymath}
             \frac{d}{dt}\iint_{\mathcal{T}} \Big[ -q\eta_{x} \Big] \,dx\,dy  -\iiint_{\mathcal{V}} \Big[ \frac{1}{2}|\nabla \phi|^{2}_{x}  + \frac{\partial}{\partial x}(-gz - \frac{1}{2}|\nabla \phi|^{2} - \frac{p}{\rho})  \Big]  dV = 0
            \end{displaymath}
            
            The \(\frac{1}{2}|\nabla \phi|^{2}_{x}\) terms cancel
            
            \begin{displaymath}
             \frac{d}{dt}\iint_{\mathcal{T}} \Big[ -q\eta_{x} \Big] \,dx\,dy  +\iiint_{\mathcal{V}} \Big[  \frac{p_{x}}{\rho} \Big]  dV = 0
            \end{displaymath}
            
            By Leibniz' rule: 
            
            \begin{displaymath}
             \frac{d}{dt}\iint_{\mathcal{T}} \Big[ -q\eta_{x} \Big] \,dx\,dy + \int_{-\infty}^{\infty}\int_{-\infty}^{\infty} \frac{d}{dx}[\int_{-h}^{\eta}
            ( \frac{p}{\rho}) dz] \,dx\,dy  -  \int_{-\infty}^{\infty}\int_{-\infty}^{\infty}  \frac{p(x,y,\eta,t)}{\rho}\eta_{x} \,dx\,dy  = 0
            \end{displaymath}
            
            Noting pressure vanishes at the surface, as well as the fact that total derivatives vanish on the surface by decay assumptions, the two right-most terms vanish
            
            \begin{displaymath}
             \frac{d}{dt}\iint_{\mathcal{T}} \Big[ -q\eta_{x} \Big] \,dx\,dy  = 0
            \end{displaymath}
            
            And we're done.

    \subsection{Conservation Law 2: \(-q\eta_{y}\)}
            
            We will omit this proof because of the symmetry of this law to law 1. Following all the same steps of law 1 but instead using the y version of the harmonic test function, noting all the variables in the integrand will switch from x to y will yield law 2.

    \subsection{Conservation Law 3: \(\mathcal{H} = \frac{1}{2} q \eta_{t} + \frac{1}{2}  g\eta^{2}\)}
            
            Here we use integral equation of motion (\ref{eqn:secondInt}) and the imaginary part of  harmonic test function \(\psi_{3}^{x}\). 
            
            \begin{displaymath}
             \frac{d}{dt}\iint_{\mathcal{T}} \Big[ \psi_{z} \eta_{t} \Big] \,dx\,dy - \iint_{ \mathcal{T}}\Big[  q_{t}  \nabla \psi_{z} +  ( \nabla^{\perp} q \cdot \nabla^{\perp} \psi_{z} ) \nabla \phi \Big] \cdot \vec{n}\, dS=  - \iint_{\mathcal{B}}\Big[\phi_{t} \psi_{zz}  \Big]  \,dx\,dy 
            \end{displaymath}

            \begin{displaymath}
             \frac{d}{dt}\iint_{\mathcal{T}} \Big[ \frac{1}{2}(x + iz)^{2}i \eta_{t} \Big] \,dx\,dy - \iint_{ \mathcal{T}}\Big[  q_{t}  \nabla \frac{1}{2}(x + iz)^{2}i +  ( \nabla^{\perp} q \cdot \nabla^{\perp} \frac{1}{2}(x + iz)^{2}i ) \nabla \phi \Big] \cdot \vec{n}\, dS=  - \iint_{\mathcal{B}}\Big[\phi_{t} (-x -iz)  \Big]  \,dx\,dy 
            \end{displaymath}
            
            Isolating imaginary part
            
            \begin{displaymath}
             \frac{d}{dt}\iint_{\mathcal{T}} \Big[ i\frac{1}{2}(x^{2} - \eta^{2})\eta_{t} \Big] \,dx\,dy - \iint_{ \mathcal{T}}\Big[  q_{t}  \nabla (i\frac{1}{2}(x^{2} - \eta^{2})) +  ( \nabla^{\perp} q \cdot \nabla^{\perp} (i\frac{1}{2}(x^{2} - \eta^{2})) \nabla \phi \Big] \cdot \vec{n}\, dS=  - \iint_{\mathcal{B}}\Big[\phi_{t} (-iz)  \Big]  \,dx\,dy 
            \end{displaymath}
            
            Dropping i's and expanding del operators
            
            \begin{displaymath}
             \frac{d}{dt}\iint_{\mathcal{T}} \Big[ \frac{1}{2}(x^{2} - \eta^{2})\eta_{t} \Big] \,dx\,dy - \iint_{ \mathcal{T}}\Big[  -q_{t}x\eta_{x} - q_{t}\eta +  q_{x}x\eta_{t} \Big] \,dx\,dy =  - \iint_{\mathcal{B}}\Big[\phi_{t}h \Big]  \,dx\,dy 
            \end{displaymath}

            \begin{displaymath}
             \frac{d}{dt}\iint_{\mathcal{T}} \Big[ \frac{1}{2}(x^{2} - \eta^{2})\eta_{t} \Big] \,dx\,dy + \iint_{ \mathcal{T}}\Big[ q_{t}\eta + (q_{t}\eta_{x} - q_{x}\eta_{t})x \Big] \,dx\,dy =  - \iint_{\mathcal{B}}\Big[\phi_{t}h \Big]  \,dx\,dy 
            \end{displaymath}
            
            Pushing to the boundary for the left-most term, using Green's second identity noting harmonic terms, then pushing back to the surface and bottom
            
            \begin{displaymath}
             \frac{d}{dt}\iint_{\mathcal{T}} \Big[ -\eta_{x}xq -q \eta \Big] \,dx\,dy  + \iint_{ \mathcal{T}}\Big[ q_{t}\eta + (q_{t}\eta_{x} - q_{x}\eta_{t})x \Big] \,dx\,dy =   0
            \end{displaymath}
            
            Note that the bottom domain terms cancel. Commuting the time derivative on the left term with the integrals assuming well-behavedness, also noting that \( -\frac{d}{dt}(q\eta) + q_{t}\eta = -q\eta_{t}\)
            
            \begin{displaymath}
            \iint_{\mathcal{T}} \Big[ -  \frac{d}{dt}(\eta_{x}xq)  -q\eta_{t} + (q_{t}\eta_{x} - q_{x}\eta_{t})x \Big] \,dx\,dy = 0
            \end{displaymath}
            
            Note that  \( -\frac{d}{dt}(\eta_{x}xq) + q_{t}\eta_{x}x = -q\eta_{xt}x\)
            
            \begin{displaymath}
            \iint_{\mathcal{T}} \Big[   -q\eta_{t} -q\eta_{xt}x - q_{x}\eta_{t}x \Big] \,dx\,dy =  0
            \end{displaymath}
            
            Adding and subtracting \(g\eta^{2}\), we recover \(-2\mathcal{H}\) on the left
            
            \begin{displaymath}
            \iint_{\mathcal{T}} \Big[   -2\mathcal{H} + g\eta^{2} -q\eta_{xt}x - q_{x}\eta_{t}x \Big] \,dx\,dy =  0
            \end{displaymath}
            
            Rearranging and dividing through by -2
            
            \begin{displaymath}
            \iint_{\mathcal{T}} \Big[ \mathcal{H} \Big] \,dx\,dy = \frac{1}{2}\iint_{\mathcal{T}} \Big[    g\eta^{2} - q\eta_{xt}x - q_{x}\eta_{t}x \Big] \,dx\,dy 
            \end{displaymath}
            
            Taking a time derivative of both sides, as well as reversing an x product rule on the right top domain integral:
            
            \begin{displaymath}
            \frac{d}{dt}\iint_{\mathcal{T}} \Big[ \mathcal{H} \Big] \,dx\,dy = \frac{d}{dt}\iint_{\mathcal{T}} \Big[  \frac{1}{2}  g\eta^{2}  - \frac{1}{2}q\eta_{xt}x -\frac{1}{2}\frac{d}{dx}(\phi \eta_{t}x) + \frac{1}{2}q\eta_{xt}x + \frac{1}{2}q\eta_{t} \Big] \,dx\,dy 
            \end{displaymath}
            
            Noting that the total derivative term vanishes, and the \(\eta_{xt}\) terms cancel
            
            \begin{equation}\label{3a}
            \frac{d}{dt}\iint_{\mathcal{T}} \Big[ \mathcal{H} \Big] \,dx\,dy = \frac{d}{dt} \frac{1}{2} \iint_{\mathcal{T}} \Big[   g\eta^{2}  +  q\eta_{t} \Big] \,dx\,dy 
            \end{equation}

            Note \(\mathcal{H} = T_{3}\). Pushing to the boundary on the RHS terms 
            
            \begin{displaymath}
            \frac{d}{dt}\iint_{\mathcal{T}} \Big[  T_{3}  \Big]  \,dx\,dy  =   \oiint_{\partial \mathcal{V}} \Big[  gz \nabla \phi \cdot \vec{n} \Big] dS  +  \frac{d}{dt}\frac{1}{2}\oiint_{\partial \mathcal{V}}\Big[  \phi \nabla \phi \cdot \vec{n} \Big] dS  
            \end{displaymath}
            
            Using the dynamic boundary condition on the first term of the RHS: 
            
            \begin{displaymath}
            \frac{d}{dt}\iint_{\mathcal{T}} \Big[  T_{3}  \Big]  \,dx\,dy  = \oiint_{\partial \mathcal{V}} \Big[ ( -\phi_{t} - \frac{1}{2}|\nabla \phi|^{2} ) \nabla \phi \cdot \vec{n} \Big] dS  +  \frac{d}{dt}\frac{1}{2}\oiint_{\partial \mathcal{V}}\Big[  \phi \nabla \phi \cdot \vec{n} \Big] dS  
            \end{displaymath}
            
            Using Green's first identity on the second term on the RHS:

            \begin{displaymath}
            \frac{d}{dt}\iint_{\mathcal{T}} \Big[  T_{3}  \Big]  \,dx\,dy  =   \oiint_{\partial \mathcal{V}} \Big[ ( -\phi_{t} - \frac{1}{2}|\nabla \phi|^{2} ) \nabla \phi \cdot \vec{n} \Big] dS  +  \frac{d}{dt}\frac{1}{2}\iiint_{\mathcal{V}}\Big[ \nabla \phi  \cdot \nabla \phi  \Big] dV  
            \end{displaymath}
            
            Using Reynold's Transport Theorem from continuum mechanics:

            \begin{displaymath}
            \frac{d}{dt}\iint_{\mathcal{T}} \Big[  T_{3}  \Big]  \,dx\,dy  =  \oiint_{\partial \mathcal{V}} \Big[ ( -\phi_{t} - \frac{1}{2}|\nabla \phi|^{2} ) \nabla \phi \cdot \vec{n} \Big] dS  + \iiint_{\mathcal{V}}\Big[ \frac{1}{2}\frac{\partial}{\partial t}\nabla \phi  \cdot \nabla \phi  \Big] dV + \oiint_{\partial \mathcal{V}} \frac{1}{2}\Big[ (\nabla \phi  \cdot \nabla \phi)\nabla \phi \cdot \vec{n} \Big] dS
            \end{displaymath}

            \begin{displaymath}
            \frac{d}{dt}\iint_{\mathcal{T}} \Big[  T_{3}  \Big]  \,dx\,dy  =  \oiint_{\partial \mathcal{V}} \Big[ ( -\phi_{t} - \frac{1}{2}|\nabla \phi|^{2} ) \nabla \phi \cdot \vec{n} \Big] dS  + \iiint_{\mathcal{V}}\Big[ \nabla \phi_{t}  \cdot \nabla \phi  \Big] dV + \oiint_{\partial \mathcal{V}} \frac{1}{2}\Big[ (\nabla \phi  \cdot \nabla \phi)\nabla \phi \cdot \vec{n} \Big] dS
            \end{displaymath}

            Pushing the bulk term back to the boundary:
            
            \begin{displaymath}
            \frac{d}{dt}\iint_{\mathcal{T}} \Big[  T_{3}  \Big]  \,dx\,dy  =   \oiint_{\partial \mathcal{V}} \Big[ ( -\phi_{t} - \frac{1}{2}|\nabla \phi|^{2} ) \nabla \phi \cdot \vec{n} \Big] dS  + \oiint_{\partial \mathcal{V}} \Big[ \phi_{t} \nabla \phi \cdot \vec{n} \Big] dS + \oiint_{\partial \mathcal{V}} \frac{1}{2}\Big[ (\nabla \phi  \cdot \nabla \phi)\nabla \phi \cdot \vec{n} \Big] dS
            \end{displaymath}
            
            The first and third term on the RHS cancel. The second and fourth also cancel, leaving us with

            \begin{displaymath}
            \frac{d}{dt}\iint_{\mathcal{T}} \Big[  T_{3}  \Big]  \,dx\,dy  =   0
            \end{displaymath}

    \subsection{Conservation Law 4: \(\eta\)}

            Substituting \(\psi_{1}^{x} = x + iz\) into (\ref{eqn:firstInt}), yields 
            \begin{displaymath}
                \frac{d}{dt} \iint_{\mathcal{T}} \Big[ (x + iz) \Big]  \,dx\,dy - \iint_{\mathcal{T}} \Big[  q_{x} \Big] \,dx\,dy =  + \iint_{\Gamma} \Big[  \phi(\nabla (i) \cdot \vec{n}) \Big] dS
            \end{displaymath}

            \begin{displaymath}
                \frac{d}{dt} \iint_{\mathcal{T}} \Big[ (x + iz) \Big]  \,dx\,dy - \iint_{\mathcal{T}} \Big[  q_{x} \Big] \,dx\,dy =  + \iint_{\Gamma} \Big[  0 \Big] dS
            \end{displaymath}

            \begin{displaymath}
                \frac{d}{dt} \iint_{\mathcal{T}} \Big[ (x + iz) \Big]  \,dx\,dy - \iint_{\mathcal{T}} \Big[  q_{x} \Big] \,dx\,dy = 0
            \end{displaymath}
            
            Isolating imaginary part:
            
            \begin{displaymath}
                \frac{d}{dt} \iint_{\mathcal{T}} \Big[ \eta \Big]  = 0
            \end{displaymath}
            
            And we're done.

    \subsection{Conservation Law 5: \(q + tg\eta\)}
            
            Starting from (\ref{eqn:firstInt}), we use \(\psi_{2}^{x}\) and isolate the real part of the resulting equation to prove this conservation law. Equation (\ref{eqn:firstInt}) reads:

            \begin{displaymath}
                \frac{d}{dt} \iint_{\mathcal{T}} \Big[ \psi \Big]  \,dx\,dy - \iint_{\mathcal{T}} \Big[  q_{x} \psi_{x}  +   q_{y}\psi_{y} \Big] \,dx\,dy =  + \iint_{\Gamma} \Big[  \phi(\nabla \psi_{z} \cdot \vec{n}) \Big] dS
            \end{displaymath}

            \begin{displaymath}
                \frac{d}{dt} \iint_{\mathcal{T}} \Big[ \psi \Big]  \,dx\,dy = \iint_{\mathcal{T}} \Big[  q_{x} \psi_{x}  +   q_{y}\psi_{y} \Big] \,dx\,dy  + \iint_{\Gamma} \Big[  \phi(\nabla \psi_{z} \cdot \vec{n}) \Big] dS
            \end{displaymath}
            
            Inserting \(\psi_{2}^{x}\)
            
            \begin{displaymath}
                \frac{d}{dt} \iint_{\mathcal{T}} \Big[ \frac{1}{2}(x^{2} - z^{2}) + izx \Big]  \,dx\,dy = \iint_{\mathcal{T}} \Big[  q_{x} (x + iz) \Big] \,dx\,dy  + \iint_{\Gamma} \Big[  \phi(\nabla \psi_{z} \cdot \vec{n}) \Big] dS
            \end{displaymath}
            
            Isolating real part:
            
            \begin{displaymath}
                \frac{d}{dt} \iint_{\mathcal{T}} \Big[ \frac{1}{2}(x^{2} - z^{2}) \Big]  \,dx\,dy = \iint_{\mathcal{T}} \Big[  q_{x}x\Big] \,dx\,dy  - \iint_{\Gamma} \Big[  \phi(\nabla z \cdot \vec{n}) \Big] dS
            \end{displaymath}
            
            The time derivative of the x term on the left-most term vanishes
            
            \begin{displaymath}
                \frac{d}{dt} \iint_{\mathcal{T}} \Big[ -\frac{1}{2}\eta^{2} \Big]  \,dx\,dy = \iint_{\mathcal{T}} \Big[  q_{x}x\Big] \,dx\,dy  - \iint_{\Gamma} \Big[  \phi(\nabla z \cdot \vec{n}) \Big] dS
            \end{displaymath}

            \begin{displaymath}
                - \iint_{\mathcal{T}} \Big[ \eta\eta_{t} \Big]  \,dx\,dy = \iint_{\mathcal{T}} \Big[  q_{x}x\Big] \,dx\,dy  - \iint_{\Gamma} \Big[  \phi(\nabla z \cdot \vec{n}) \Big] dS
            \end{displaymath}

            \begin{displaymath}
                - \oiint_{\partial \mathcal{V}} \Big[ z \nabla \phi  \cdot \vec{n} \Big]  dS = \iint_{\mathcal{T}} \Big[  q_{x}x \Big] \,dx\,dy  - \iint_{\Gamma} \Big[  \phi(\nabla z \cdot \vec{n}) \Big] dS
            \end{displaymath}

            \begin{displaymath}
              - \iint_{\mathcal{T}} \Big[ z \nabla \phi  \cdot \vec{n} \Big]  dS   - \iint_{\Gamma} \Big[ z \nabla \phi  \cdot \vec{n} \Big]  dS = \iint_{\mathcal{T}} \Big[  q_{x}x \Big] \,dx\,dy  - \iint_{\Gamma} \Big[  \phi(\nabla z \cdot \vec{n}) \Big] dS
            \end{displaymath}
            
            The \(\Gamma\) terms cancel by Green's Second Identity:
            
            \begin{displaymath}
              - \iint_{\mathcal{T}} \Big[ z \nabla \phi  \cdot \vec{n} \Big]  dS   = \iint_{\mathcal{T}} \Big[  q_{x}x \Big] \,dx\,dy  
            \end{displaymath}
            
            But due to the bed condition and the outer normal on the bottom surface causing the bottom contribution to vanish, we can push the LHS to the boundary:
            
            \begin{displaymath}
              - \oiint_{\partial \mathcal{V}} \Big[ z \nabla \phi  \cdot \vec{n} \Big]  dS   = \iint_{\mathcal{T}} \Big[  q_{x}x \Big] \,dx\,dy  
            \end{displaymath}
            
            By Green's first identity on the LHS: 
            
            \begin{displaymath}
              - \iiint_{\mathcal{V}} \Big[ \nabla z \cdot \nabla \phi  \Big] dV  = \iint_{\mathcal{T}} \Big[  q_{x}x \Big] \,dx\,dy  
            \end{displaymath}
            
            Reversing a product rule on the RHS:
            
            \begin{displaymath}
              - \iiint_{\mathcal{V}} \Big[ \nabla z \cdot \nabla \phi  \Big] dV  = \iint_{\mathcal{T}} \Big[  \frac{d}{dx}(\phi x) - q \Big] \,dx\,dy  
            \end{displaymath}
            
            The total derivative term vanishes. Taking a time derivative of both sides
            
            \begin{displaymath}
              - \frac{d}{dt} \iiint_{\mathcal{V}} \Big[ \frac{\partial}{\partial z} \phi \Big] dV  = \frac{d}{dt}\iint_{\mathcal{T}} \Big[   - q \Big] \,dx\,dy  
            \end{displaymath}
            
            By Reynold's transport theorem from continuum mechanics:
            
            \begin{displaymath}
              - \iiint_{\mathcal{V}} \Big[ \frac{\partial}{\partial z} \phi_{t} \Big] dV - \oiint_{\partial \mathcal{V}} (\nabla \phi \cdot \nabla z) \nabla \phi \cdot \vec{n}\, dS  = -\frac{d}{dt}\iint_{\mathcal{T}} \Big[  q \Big] \,dx\,dy  
            \end{displaymath}
            
            The dynamic boundary condition on the LHS gives us 
            
            \begin{displaymath}
               \iiint_{\mathcal{V}} \Big[ \frac{\partial}{\partial z} (g\eta + \frac{1}{2}|\nabla \phi|^{2}) \Big] dV - \oiint_{\partial \mathcal{V}} \phi_{z} \nabla \phi \cdot \vec{n}\, dS  = -\frac{d}{dt}\iint_{\mathcal{T}} \Big[  q \Big] \,dx\,dy  
            \end{displaymath}
            
            Green's first identity on the second term on the LHS gives us:
            
            \begin{displaymath}
               \iiint_{\mathcal{V}} \Big[ \frac{\partial}{\partial z} (g\eta + \frac{1}{2}|\nabla \phi|^{2}) \Big] dV - \iiint_{\mathcal{V}} (\nabla \phi_{z} \cdot \nabla \phi) dV  = -\frac{d}{dt}\iint_{\mathcal{T}} \Big[  q \Big] \,dx\,dy  
            \end{displaymath}
            
            Bringing the \(  \frac{1}{2}|\nabla \phi|^{2} \) term on the left-most integral to the next integral on the right
            
            \begin{displaymath}
               \iiint_{\mathcal{V}} \Big[ \frac{\partial}{\partial z} g\eta  \Big] dV - \iiint_{\mathcal{V}} (\nabla \phi_{z} \cdot \nabla \phi - \frac{\partial}{\partial z}\frac{1}{2}|\nabla \phi|^{2}) dV  = -\frac{d}{dt}\iint_{\mathcal{T}} \Big[  q \Big] \,dx\,dy  
            \end{displaymath}

            \begin{displaymath}
               \int_{-\infty}^{\infty}\int_{-\infty}^{\infty} \int_{-h}^{\eta}
             \Big[  g  \Big] dzdydx - \iiint_{\mathcal{V}} (\nabla \phi_{z} \cdot \nabla \phi - \frac{\partial}{\partial z}\frac{1}{2}|\nabla \phi|^{2}) dV  = -\frac{d}{dt}\iint_{\mathcal{T}} \Big[  q \Big] \,dx\,dy  
            \end{displaymath}

            \begin{displaymath}
               \int_{-\infty}^{\infty}\int_{-\infty}^{\infty}
             \Big[  g\eta + gh  \Big] dydx - \iiint_{\mathcal{V}} (\nabla \phi_{z} \cdot \nabla \phi - \frac{\partial}{\partial z}\frac{1}{2}|\nabla \phi|^{2}) dV  = -\frac{d}{dt}\iint_{\mathcal{T}} \Big[  q \Big] \,dx\,dy  
            \end{displaymath}
            
            Splitting up the first term on the LHS into top and bottom contributions and then bringing the bottom contribution to the RHS
            
            \begin{displaymath}
               \iint_{\mathcal{T}} \Big[  g\eta \Big] dydx - \iiint_{\mathcal{V}} (\nabla \phi_{z} \cdot \nabla \phi - \frac{\partial}{\partial z}\frac{1}{2}|\nabla \phi|^{2}) dV  = -\frac{d}{dt}\iint_{\mathcal{T}} \Big[  q \Big] \,dx\,dy  -  \iint_{\mathcal{B}}
             \Big[  gh \Big] dydx
            \end{displaymath}
            
            Expanding the second term on the LHS
            
            \begin{displaymath}
               \iint_{\mathcal{T}} \Big[  g\eta \Big] dydx - \iiint_{\mathcal{V}} (\nabla \phi_{z} \cdot \nabla \phi - \frac{\partial}{\partial z}\frac{1}{2}\nabla \phi \cdot \nabla \phi) dV  = -\frac{d}{dt}\iint_{\mathcal{T}} \Big[  q \Big] \,dx\,dy  -  \iint_{\mathcal{B}}
             \Big[  gh \Big] dydx
            \end{displaymath}

            \begin{displaymath}
               \iint_{\mathcal{T}} \Big[  g\eta \Big] dydx - \iiint_{\mathcal{V}} (\nabla \phi_{z} \cdot \nabla \phi - \frac{1}{2}\nabla \phi_{z} \cdot \nabla \phi  - \frac{1}{2}\nabla \phi \cdot \nabla \phi_{z})dV  = -\frac{d}{dt}\iint_{\mathcal{T}} \Big[  q \Big] \,dx\,dy  -  \iint_{\mathcal{B}}
             \Big[  gh \Big] dydx
            \end{displaymath}

            \begin{displaymath}
               \iint_{\mathcal{T}} \Big[  g\eta \Big] dydx - \iiint_{\mathcal{V}} (\nabla \phi_{z} \cdot \nabla \phi - \nabla \phi_{z} \cdot \nabla \phi)dV  = -\frac{d}{dt}\iint_{\mathcal{T}} \Big[  q \Big] \,dx\,dy  -  \iint_{\mathcal{B}}
             \Big[  gh \Big] dydx
            \end{displaymath}

            \begin{displaymath}
               \iint_{\mathcal{T}} \Big[  g\eta \Big] dydx - \iiint_{\mathcal{V}} 0 dV  = -\frac{d}{dt}\iint_{\mathcal{T}} \Big[  q \Big] \,dx\,dy  -  \iint_{\mathcal{B}}
             \Big[  gh \Big] dydx
            \end{displaymath}
            
            Re-arranging and reversing a time derivative product rule
            
            \begin{displaymath}
              \frac{d}{dt}\iint_{\mathcal{T}} \Big[  q \Big] \,dx\,dy +  \iint_{\mathcal{T}} \Big[  \frac{d}{dt}(tg\eta) - tg\eta_{t} \Big] dydx   =   -  \iint_{\mathcal{B}}
             \Big[  gh \Big] dydx
            \end{displaymath}

            \begin{displaymath}
              \frac{d}{dt}\iint_{\mathcal{T}} \Big[  q \Big] \,dx\,dy +  \iint_{\mathcal{T}} \Big[  \frac{d}{dt}(tg\eta) - tg  \frac{d}{dt}\iint_{\mathcal{T}} \Big[  \eta \Big] dydx   =   -  \iint_{\mathcal{B}}
             \Big[  gh \Big] dydx
            \end{displaymath}
            
            Recognizing law 4 on the LHS, it vanishes
            
            \begin{displaymath}
              \frac{d}{dt}\iint_{\mathcal{T}} \Big[  q \Big] \,dx\,dy +  \iint_{\mathcal{T}} \Big[  \frac{d}{dt}(tg\eta) \Big] \,dx\,dy - tg \cancelto{0}{ \frac{d}{dt}\iint_{\mathcal{T}} \Big[  T_{4} \Big] dydx  } =   -  \iint_{\mathcal{B}}
             \Big[  gh \Big] dydx
            \end{displaymath}

            \begin{displaymath}
              \frac{d}{dt}\iint_{\mathcal{T}} \Big[  q  + tg\eta \Big] \,dx\,dy =   -  \iint_{\mathcal{B}}
             \Big[  gh \Big] dydx
            \end{displaymath}
        
    \subsection{Conservation Law 6: \(x\eta + tq\eta_{x}\)}
        
        % Starting from (\ref{eqn:firstInt}), we use \(\psi_{2}^{x}\) and isolate the imaginary part of the resulting equation to prove this conservation law. Equation (\ref{eqn:firstInt}) reads:
         
         \begin{displaymath}
             \frac{d}{dt} \iint_{\mathcal{T}} \Big[ \psi \Big]  \,dx\,dy - \iint_{\mathcal{T}} \Big[  q_{x} \psi_{x}  +   q_{y}\psi_{y} \Big] \,dx\,dy =  + \iint_{\Gamma} \Big[  \phi(\nabla \psi_{z} \cdot \vec{n}) \Big] dS
         \end{displaymath}
         
         Plugging in \(\psi_{2}^{x}\) 
         
         \begin{displaymath}
             \frac{d}{dt} \iint_{\mathcal{T}} \Big[ \frac{1}{2}(x + iz)^{2} \Big]  \,dx\,dy - \iint_{\mathcal{T}} \Big[  q_{x} (x + iz)  \Big] \,dx\,dy =  + \iint_{\Gamma} \Big[  \phi(\nabla (-z + ix) \cdot \vec{n}) \Big] dS
         \end{displaymath}
         
         Expanding the LHS
         
         \begin{displaymath}
             \frac{d}{dt} \iint_{\mathcal{T}} \Big[ \frac{1}{2}(x^{2} - {z}^{2}) + ixz \Big]  \,dx\,dy - \iint_{\mathcal{T}} \Big[  q_{x} (x + iz)  \Big] \,dx\,dy =  + \iint_{\Gamma} \Big[  \phi(\nabla (-z + ix) \cdot \vec{n}) \Big] dS
         \end{displaymath}
         
         Isolating imaginary part and dropping i's:

         \begin{displaymath}
             \frac{d}{dt} \iint_{\mathcal{T}} \Big[ xz \Big]  \,dx\,dy - \iint_{\mathcal{T}} \Big[  q_{x}\eta   \Big] \,dx\,dy =  + \iint_{\Gamma} \Big[  \phi(\nabla (x) \cdot \vec{n}) \Big] dS
         \end{displaymath}
         
         We previously argued that the term on the RHS vanishes in law 1:

         \begin{displaymath}
             \frac{d}{dt} \iint_{\mathcal{T}} \Big[ xz \Big]  \,dx\,dy - \iint_{\mathcal{T}} \Big[  q_{x}\eta  \Big] \,dx\,dy =  0
         \end{displaymath}

         \begin{displaymath}
             \frac{d}{dt} \iint_{\mathcal{T}} \Big[ xz \Big]  \,dx\,dy - \iint_{\mathcal{T}} \Big[  \frac{d}{dx}(\phi \eta) - q\eta_{x}  \Big] \,dx\,dy =  0
         \end{displaymath}
         
         On the RHS, the total derivative term vanishes:
         
         \begin{displaymath}
             \frac{d}{dt} \iint_{\mathcal{T}} \Big[ xz \Big]  \,dx\,dy - \iint_{\mathcal{T}} \Big[  - q\eta_{x}  \Big] \,dx\,dy =  0
         \end{displaymath}
         
         Inverting a time derivative:
         
         \begin{displaymath}
             \frac{d}{dt} \iint_{\mathcal{T}} \Big[ xz \Big]  \,dx\,dy  \iint_{\mathcal{T}} \Big[   \frac{d}{dt}(tq\eta_{x}) - t\frac{d}{dt}(q\eta_{x})  \Big] \,dx\,dy =  0
         \end{displaymath}

         \begin{displaymath}
             \frac{d}{dt} \iint_{\mathcal{T}} \Big[ xz \Big]  \,dx\,dy + \iint_{\mathcal{T}} \Big[   \frac{d}{dt}(tq\eta_{x}) \Big] \,dx\,dy - t \iint_{\mathcal{T}} \Big[ \frac{d}{dt}(q\eta_{x})  \Big] \,dx\,dy =  0
         \end{displaymath}

         \begin{displaymath}
             \frac{d}{dt} \iint_{\mathcal{T}} \Big[ xz \Big]  \,dx\,dy + \iint_{\mathcal{T}} \Big[   \frac{d}{dt}(tq\eta_{x}) \Big] \,dx\,dy + t \cancelto{0}{\frac{d}{dt} \iint_{\mathcal{T}} \Big[ T_{1}  \Big] \,dx\,dy} =  0
         \end{displaymath}

         \begin{displaymath}
             \frac{d}{dt} \iint_{\mathcal{T}} \Big[ x\eta \Big]  \,dx\,dy +  \frac{d}{dt}  \iint_{\mathcal{T}} \Big[   tq\eta_{x} \Big] \,dx\,dy  =  0
         \end{displaymath}

         \begin{displaymath}
             \frac{d}{dt} \iint_{\mathcal{T}} [ x\eta +  tq\eta_{x} ] \,dx\,dy = 0
         \end{displaymath}
         
         As required.

    \subsection{Conservation Law 7: \(y\eta + tq\eta_{y}\)}
         
         Starting from (\ref{eqn:firstInt}), we use \(\psi_{2}^{y}\) and isolate the imaginary part of the resulting equation to prove this conservation law. Equation of motion (\ref{eqn:firstInt}) reads:
         
         \begin{displaymath}
             \frac{d}{dt} \iint_{\mathcal{T}} \Big[ \psi \Big]  \,dx\,dy - \iint_{\mathcal{T}} \Big[  q_{x} \psi_{x}  +   q_{y}\psi_{y} \Big] \,dx\,dy =  + \iint_{\Gamma} \Big[  \phi(\nabla \psi_{z} \cdot \vec{n}) \Big] dS
         \end{displaymath}
         
         Plugging in \(\psi_{2}^{y}\) 
         
         \begin{displaymath}
             \frac{d}{dt} \iint_{\mathcal{T}} \Big[ \frac{1}{2}(y + iz)^{2} \Big]  \,dx\,dy - \iint_{\mathcal{T}} \Big[  q_{y} (y + iz)  \Big] \,dx\,dy =  + \iint_{\Gamma} \Big[  \phi(\nabla (-z + iy) \cdot \vec{n}) \Big] dS
         \end{displaymath}
         
         Expanding the LHS
         
         \begin{displaymath}
             \frac{d}{dt} \iint_{\mathcal{T}} \Big[ \frac{1}{2}(y^{2} - {z}^{2}) + iyz \Big]  \,dx\,dy - \iint_{\mathcal{T}} \Big[  q_{y} (y + iz)  \Big] \,dx\,dy =  + \iint_{\Gamma} \Big[  \phi(\nabla (-z + iy) \cdot \vec{n}) \Big] dS
         \end{displaymath}
         
         Isolating imaginary part and dropping i's:

         \begin{displaymath}
             \frac{d}{dt} \iint_{\mathcal{T}} \Big[ yz \Big]  \,dx\,dy - \iint_{\mathcal{T}} \Big[  q_{y}\eta   \Big] \,dx\,dy =  + \iint_{\Gamma} \Big[  \phi(\nabla (y) \cdot \vec{n}) \Big] dS
         \end{displaymath}
         
         We previously argued that the term on the RHS vanishes in law 2:

         \begin{displaymath}
             \frac{d}{dt} \iint_{\mathcal{T}} \Big[ yz \Big]  \,dx\,dy - \iint_{\mathcal{T}} \Big[  q_{y}\eta  \Big] \,dx\,dy =  0
         \end{displaymath}

         \begin{displaymath}
             \frac{d}{dt} \iint_{\mathcal{T}} \Big[ yz \Big]  \,dx\,dy - \iint_{\mathcal{T}} \Big[  \frac{d}{dy}(\phi \eta) - q\eta_{y}  \Big] \,dx\,dy =  0
         \end{displaymath}
         
         On the RHS, the total derivative term vanishes:
         
         \begin{displaymath}
             \frac{d}{dt} \iint_{\mathcal{T}} \Big[ yz \Big]  \,dx\,dy - \iint_{\mathcal{T}} \Big[  - q\eta_{y}  \Big] \,dx\,dy =  0
         \end{displaymath}
         
         Inverting a time derivative:
         
         \begin{displaymath}
             \frac{d}{dt} \iint_{\mathcal{T}} \Big[ yz \Big]  \,dx\,dy  \iint_{\mathcal{T}} \Big[   \frac{d}{dt}(tq\eta_{y}) - t\frac{d}{dt}(q\eta_{y})  \Big] \,dx\,dy =  0
         \end{displaymath}

         \begin{displaymath}
             \frac{d}{dt} \iint_{\mathcal{T}} \Big[ yz \Big]  \,dx\,dy + \iint_{\mathcal{T}} \Big[   \frac{d}{dt}(tq\eta_{y}) \Big] \,dx\,dy - t \iint_{\mathcal{T}} \Big[ \frac{d}{dt}(q\eta_{y})  \Big] \,dx\,dy =  0
         \end{displaymath}

         \begin{displaymath}
             \frac{d}{dt} \iint_{\mathcal{T}} \Big[ yz \Big]  \,dx\,dy + \iint_{\mathcal{T}} \Big[   \frac{d}{dt}(tq\eta_{y}) \Big] \,dx\,dy + t \cancelto{0}{\frac{d}{dt} \iint_{\mathcal{T}} \Big[ T_{2}  \Big] \,dx\,dy} =  0
         \end{displaymath}

         \begin{displaymath}
             \frac{d}{dt} \iint_{\mathcal{T}} \Big[ y\eta \Big]  \,dx\,dy +  \frac{d}{dt}  \iint_{\mathcal{T}} \Big[   tq\eta_{y} \Big] \,dx\,dy  =  0
         \end{displaymath}

         \begin{displaymath}
             \frac{d}{dt} \iint_{\mathcal{T}} [ y\eta +  tq\eta_{y} ] \,dx\,dy = 0
         \end{displaymath}
         
         As required.

    \subsection{Conservation Law 8: \(\frac{1}{2} \eta^{2} - t T_{5} + \frac{1}{2} g t^{2}  T_{4}\) }
         
         Starting from (\ref{eqn:firstInt}), we use \(\psi_{2}^{x}\) and isolate the real part of the resulting equation to prove this conservation law. Equation of motion (\ref{eqn:firstInt}) reads:

         \begin{displaymath}
             \frac{d}{dt} \iint_{\mathcal{T}} \Big[ \psi \Big]  \,dx\,dy - \iint_{\mathcal{T}} \Big[  q_{x} \psi_{x}  +   q_{y}\psi_{y} \Big] \,dx\,dy =  + \iint_{\Gamma} \Big[  \phi(\nabla \psi_{z} \cdot \vec{n}) \Big] dS
         \end{displaymath}

         Plugging in \(\psi_{2}^{x}\)
         
         \begin{displaymath}
             \frac{d}{dt} \iint_{\mathcal{T}} \Big[ \frac{1}{2}(x^{2} - z^{2}) + ixz\Big]  \,dx\,dy - \iint_{\mathcal{T}} \Big[  q_{x} (x + iz) \Big] \,dx\,dy =  + \iint_{\Gamma} \Big[  \phi(\nabla (-z + ix) \cdot \vec{n}) \Big] dS
         \end{displaymath}

         Isolating real part:

         \begin{displaymath}
             \frac{d}{dt} \iint_{\mathcal{T}} \Big[ \frac{1}{2}(x^{2} - z^{2}) \Big]  \,dx\,dy - \iint_{\mathcal{T}} \Big[  q_{x} x \Big] \,dx\,dy =  - \iint_{\Gamma} \Big[  \phi(\nabla z \cdot \vec{n}) \Big] dS
         \end{displaymath}
         
         Simplifying and reducing the RHS to the bottom surface contribution via decay arguments

         \begin{displaymath}
             \frac{d}{dt} \iint_{\mathcal{T}} \Big[ - \frac{1}{2} \eta^{2} \Big]  \,dx\,dy - \iint_{\mathcal{T}} \Big[  q_{x} x \Big] \,dx\,dy =   \iint_{\mathcal{B}} \Big[  \phi \Big] \,dx\,dy
         \end{displaymath}

         \begin{displaymath}
             \frac{d}{dt} \iint_{\mathcal{T}} \Big[ \frac{1}{2} \eta^{2} \Big]  \,dx\,dy + \iint_{\mathcal{T}} \Big[  q_{x} x \Big] \,dx\,dy =   \iint_{\mathcal{B}} \Big[ - \phi \Big] \,dx\,dy
         \end{displaymath}
         
         Reversing a product rule
         
         \begin{displaymath}
             \frac{d}{dt} \iint_{\mathcal{T}} \Big[ \frac{1}{2} \eta^{2} \Big]  \,dx\,dy + \iint_{\mathcal{T}} \Big[ \frac{d}{dx}( \phi x ) - q \Big] \,dx\,dy =   \iint_{\mathcal{B}} \Big[ - \phi \Big] \,dx\,dy
         \end{displaymath}
         
         On the RHS, the total derivative term vanishes:
         
         \begin{displaymath}
             \frac{d}{dt} \iint_{\mathcal{T}} \Big[ \frac{1}{2} \eta^{2} \Big]  \,dx\,dy - \iint_{\mathcal{T}} \Big[  q \Big] \,dx\,dy =   \iint_{\mathcal{B}} \Big[ - \phi \Big] \,dx\,dy
         \end{displaymath}
         
         Reversing a product rule
         
         \begin{displaymath}
             \frac{d}{dt} \iint_{\mathcal{T}} \Big[ \frac{1}{2} \eta^{2} \Big]  \,dx\,dy - \iint_{\mathcal{T}} \Big[  \frac{d}{dt} (qt) - tq_{t} \Big] \,dx\,dy =   \iint_{\mathcal{B}} \Big[ - \phi \Big] \,dx\,dy
         \end{displaymath}

         \begin{displaymath}
             \frac{d}{dt} \iint_{\mathcal{T}} \Big[ \frac{1}{2} \eta^{2} - qt \Big]  \,dx\,dy +   \iint_{\mathcal{T}} \Big[   tq_{t} \Big] \,dx\,dy =   \iint_{\mathcal{B}} \Big[ - \phi \Big] \,dx\,dy
         \end{displaymath}
         
         Expanding the total derivative \(q_{t} = \frac{d}{dt} \phi(x,y,\eta(x,y,t),t) = \phi_{t} + \phi_{z}\eta_{t}\): 
         
         \begin{displaymath}
             \frac{d}{dt} \iint_{\mathcal{T}} \Big[ \frac{1}{2} \eta^{2} - qt \Big]  \,dx\,dy +   t\iint_{\mathcal{T}} \Big[   \phi_{t} + \phi_{z}\eta_{t}  \Big] \,dx\,dy =   \iint_{\mathcal{B}} \Big[ - \phi \Big] \,dx\,dy
         \end{displaymath}
         
         Inserting the dynamic boundary condition:
         
         \begin{displaymath}
             \frac{d}{dt} \iint_{\mathcal{T}} \Big[ \frac{1}{2} \eta^{2} - qt \Big]  \,dx\,dy +   t\iint_{\mathcal{T}} \Big[  -g\eta - \frac{1}{2}|\nabla \phi|^{2} + \phi_{z}\eta_{t}  \Big] \,dx\,dy =   \iint_{\mathcal{B}} \Big[ - \phi \Big] \,dx\,dy
         \end{displaymath}
         
         Reversing a product rule
         
         \begin{displaymath}
             \frac{d}{dt} \iint_{\mathcal{T}} \Big[ \frac{1}{2} \eta^{2} - qt \Big]  \,dx\,dy +   \iint_{\mathcal{T}} \Big[  -\frac{d}{dt}(\frac{1}{2}t^{2}g\eta) +\frac{1}{2}t^{2}g\eta_{t}  - t\frac{1}{2}|\nabla \phi|^{2} + t\phi_{z}\eta_{t}  \Big] \,dx\,dy =   \iint_{\mathcal{B}} \Big[ - \phi \Big] \,dx\,dy
         \end{displaymath}

         \begin{displaymath}
             \frac{d}{dt} \iint_{\mathcal{T}} \Big[ \frac{1}{2} \eta^{2} - qt - \frac{1}{2}t^{2}g\eta \Big]  \,dx\,dy    +\frac{1}{2}t^{2}g \iint_{\mathcal{T}} \Big[ \eta_{t} \Big] \,dx\,dy + \iint_{\mathcal{T}} \Big[ - t\frac{1}{2}|\nabla \phi|^{2} + t\phi_{z}\eta_{t}  \Big] \,dx\,dy =   \iint_{\mathcal{B}} \Big[ - \phi \Big] \,dx\,dy
         \end{displaymath}

         \begin{displaymath}
             \frac{d}{dt} \iint_{\mathcal{T}} \Big[ \frac{1}{2} \eta^{2} - qt - \frac{1}{2}t^{2}g\eta \Big]  \,dx\,dy    +\frac{1}{2}t^{2}g \frac{d}{dt} \iint_{\mathcal{T}} \Big[ T_{4} \Big] \,dx\,dy + \iint_{\mathcal{T}} \Big[ - t\frac{1}{2}|\nabla \phi|^{2} + t\phi_{z}\eta_{t}  \Big] \,dx\,dy =   \iint_{\mathcal{B}} \Big[ - \phi \Big] \,dx\,dy
         \end{displaymath}

         Recognizing law 4 on the LHS, that term vanishes:
         
         \begin{displaymath}
             \frac{d}{dt} \iint_{\mathcal{T}} \Big[ \frac{1}{2} \eta^{2} - qt - \frac{1}{2}t^{2}g\eta \Big]  \,dx\,dy + t\iint_{\mathcal{T}} \Big[ - \frac{1}{2}|\nabla \phi|^{2} + \phi_{z}\eta_{t}  \Big] \,dx\,dy =   \iint_{\mathcal{B}} \Big[ - \phi \Big] \,dx\,dy
         \end{displaymath}
         
         Rewriting the terms under the time derivative on the LHS: 
         
         \begin{displaymath}
             \frac{d}{dt} \iint_{\mathcal{T}} \Big[ \frac{1}{2} \eta^{2} - qt - t^{2}g\eta + \frac{1}{2}t^{2}g\eta \Big]  \,dx\,dy + t\iint_{\mathcal{T}} \Big[ - \frac{1}{2}|\nabla \phi|^{2} + \phi_{z}\eta_{t}  \Big] \,dx\,dy =   \iint_{\mathcal{B}} \Big[ - \phi \Big] \,dx\,dy
         \end{displaymath}

         \begin{displaymath}
             \frac{d}{dt} \iint_{\mathcal{T}} \Big[ \frac{1}{2} \eta^{2} - t(q + tg\eta) + \frac{1}{2}t^{2}g\eta \Big]  \,dx\,dy + t\iint_{\mathcal{T}} \Big[ - \frac{1}{2}|\nabla \phi|^{2} + \phi_{z}\eta_{t}  \Big] \,dx\,dy =   \iint_{\mathcal{B}} \Big[ - \phi \Big] \,dx\,dy
         \end{displaymath}

         \begin{displaymath}
             \frac{d}{dt} \iint_{\mathcal{T}} \Big[ \frac{1}{2} \eta^{2} - tT_{5} + \frac{1}{2}t^{2}gT_{4} \Big]  \,dx\,dy + t\iint_{\mathcal{T}} \Big[ - \frac{1}{2}|\nabla \phi|^{2} + \phi_{z}\eta_{t}  \Big] \,dx\,dy =   \iint_{\mathcal{B}} \Big[ - \phi \Big] \,dx\,dy
         \end{displaymath}
         
         Pushing the remaining terms on the LHS to the bulk and applying divergence theorem:
         
         \begin{displaymath}
             \frac{d}{dt} \iint_{\mathcal{T}} \Big[ \frac{1}{2} \eta^{2} - tT_{5} + \frac{1}{2}t^{2}gT_{4} \Big]  \,dx\,dy + t\iiint_{\mathcal{V}} \nabla \cdot \Big[ \begin{pmatrix} \phi_{z}\phi_{x} \\ \phi_{z}\phi_{y}  \\ - \frac{1}{2}|\nabla \phi|^{2} + \phi_{z}^{2} \end{pmatrix}   \Big] dV - t\iint_{\Gamma} \Big[ \begin{pmatrix} \phi_{z}\phi_{x} \\ \phi_{z}\phi_{y}  \\ - \frac{1}{2}|\nabla \phi|^{2} + \phi_{z}^{2} \end{pmatrix}  \Big] \cdot \vec{n} \,dS =   \iint_{\mathcal{B}} \Big[ - \phi \Big] \,dx\,dy
         \end{displaymath}

         \begin{displaymath}
             \frac{d}{dt} \iint_{\mathcal{T}} \Big[ \frac{1}{2} \eta^{2} - tT_{5} + \frac{1}{2}t^{2}gT_{4} \Big]  \,dx\,dy + t\iiint_{\mathcal{V}} \nabla \cdot \Big[ \begin{pmatrix} \phi_{z}\phi_{x} \\ \phi_{z}\phi_{y}  \\ - \frac{1}{2}|\nabla \phi|^{2} + \phi_{z}^{2} \end{pmatrix}   \Big] dV = t\iint_{\Gamma} \Big[ \begin{pmatrix} \phi_{z}\phi_{x} \\ \phi_{z}\phi_{y}  \\ - \frac{1}{2}|\nabla \phi|^{2} + \phi_{z}^{2} \end{pmatrix}  \Big] \cdot \vec{n} \,dS  -  \iint_{\mathcal{B}} \Big[ \phi \Big] \,dx\,dy
         \end{displaymath}
         
         Reducing the RHS to the bottom surface via decay assumptions:
         
         \begin{displaymath}
             \frac{d}{dt} \iint_{\mathcal{T}} \Big[ \frac{1}{2} \eta^{2} - tT_{5} + \frac{1}{2}t^{2}gT_{4} \Big]  \,dx\,dy + t\iiint_{\mathcal{V}} \nabla \cdot \Big[ \begin{pmatrix} \phi_{z}\phi_{x} \\ \phi_{z}\phi_{y}  \\ - \frac{1}{2}|\nabla \phi|^{2} + \phi_{z}^{2} \end{pmatrix}   \Big] dV = \iint_{\mathcal{B}}  t( \frac{1}{2}|\nabla \phi|^{2} - \phi_{z}^{2}) - \phi  \Big]  \,dx\,dy 
         \end{displaymath}
         
         Expanding the divergence on the LHS noting harmonicity of \(\phi\) causes terms to vanish

         \begin{displaymath}
             \frac{d}{dt} \iint_{\mathcal{T}} \Big[ \frac{1}{2} \eta^{2} - tT_{5} + \frac{1}{2}t^{2}gT_{4} \Big]  \,dx\,dy + t\iiint_{\mathcal{V}} \Big[  \phi_{zx}\phi_{x} + \phi_{zy}\phi_{y}  - \frac{\partial}{\partial z}(\frac{1}{2}|\nabla \phi|^{2}) + \phi_{zz}\phi_{z} \Big] dV = \iint_{\mathcal{B}}  t( \frac{1}{2}|\nabla \phi|^{2} - \phi_{z}^{2}) - \phi  \Big]  \,dx\,dy 
         \end{displaymath}
         
         Re-arranging
         
         \begin{displaymath}
             \frac{d}{dt} \iint_{\mathcal{T}} \Big[ \frac{1}{2} \eta^{2} - tT_{5} + \frac{1}{2}t^{2}gT_{4} \Big]  \,dx\,dy + t\iiint_{\mathcal{V}} \Big[  \phi_{zx}\phi_{x} + \phi_{zy}\phi_{y} + \phi_{zz}\phi_{z} - \frac{\partial}{\partial z}(\frac{1}{2}|\nabla \phi|^{2})  \Big] dV = \iint_{\mathcal{B}}  t( \frac{1}{2}|\nabla \phi|^{2} - \phi_{z}^{2}) - \phi  \Big]  \,dx\,dy 
         \end{displaymath}

         \begin{displaymath}
             \frac{d}{dt} \iint_{\mathcal{T}} \Big[ \frac{1}{2} \eta^{2} - tT_{5} + \frac{1}{2}t^{2}gT_{4} \Big]  \,dx\,dy + t\iiint_{\mathcal{V}} \Big[  \frac{\partial}{\partial z}( \frac{1}{2}\phi_{x}^{2} + \frac{1}{2}\phi_{y}^{2} +  \frac{1}{2}\phi_{z}^{2} ) - \frac{\partial}{\partial z}(\frac{1}{2}|\nabla \phi|^{2})  \Big] dV = \iint_{\mathcal{B}}  t( \frac{1}{2}|\nabla \phi|^{2} - \phi_{z}^{2}) - \phi  \Big]  \,dx\,dy 
         \end{displaymath}

         \begin{displaymath}
             \frac{d}{dt} \iint_{\mathcal{T}} \Big[ \frac{1}{2} \eta^{2} - tT_{5} + \frac{1}{2}t^{2}gT_{4} \Big]  \,dx\,dy + t\iiint_{\mathcal{V}} \Big[ \frac{\partial}{\partial z}( \frac{1}{2}|\nabla \phi|^{2}  ) - \frac{\partial}{\partial z}(\frac{1}{2}|\nabla \phi|^{2})  \Big] dV = \iint_{\mathcal{B}}  t( \frac{1}{2}|\nabla \phi|^{2} - \phi_{z}^{2}) - \phi  \Big]  \,dx\,dy 
         \end{displaymath}

         \begin{displaymath}
             \frac{d}{dt} \iint_{\mathcal{T}} \Big[ \frac{1}{2} \eta^{2} - tT_{5} + \frac{1}{2}t^{2}gT_{4} \Big]  \,dx\,dy + t\iiint_{\mathcal{V}} \Big[  0 \Big] dV = \iint_{\mathcal{B}}  t( \frac{1}{2}|\nabla \phi|^{2} - \phi_{z}^{2}) - \phi  \Big]  \,dx\,dy 
         \end{displaymath}

         \begin{displaymath}
             \frac{d}{dt} \iint_{\mathcal{T}} \Big[ \frac{1}{2} \eta^{2} - tT_{5} + \frac{1}{2}t^{2}gT_{4} \Big]  \,dx\,dy  = \iint_{\mathcal{B}}  t( \frac{1}{2}|\nabla \phi|^{2} - \phi_{z}^{2}) - \phi  \Big]  \,dx\,dy 
         \end{displaymath}
         
         By the bed condition:
         
         \begin{displaymath}
             \frac{d}{dt} \iint_{\mathcal{T}} \Big[ \frac{1}{2} \eta^{2} - tT_{5} + \frac{1}{2}t^{2}gT_{4} \Big]  \,dx\,dy  = \iint_{\mathcal{B}}  t \frac{1}{2}|\nabla^{\perp} \phi|^{2} - \phi  \Big]  \,dx\,dy 
         \end{displaymath}
         
         And we're done.

    \subsection{Conservation Law 9: \(x\eta_{y} - y\eta_{x})q\)}
         
             This conservation law comes out of isolating the imaginary part of integral equation of motion (\ref{eqn:firstInt}) using harmonic test functions \(\psi_{2}^{xy}(x,y)\). Starting from integral equation of motion (\ref{eqn:firstInt}) and integrating both sides by z:
         
             \begin{displaymath}
                 \frac{d}{dt} \iint_{\mathcal{T}} \int \Big[\frac{(x + iy)^2}{2}\Big] dz \,dx\,dy - \iint_{\mathcal{T}} \int \Big[  q_{x} (x + iy)  +   q_{y}(-y + ix) \Big] \,dz \,dx\,dy =  + \iint_{\Gamma} \int  \Big[  \phi (\nabla \psi_{z} \cdot \vec{n}) \Big] \,dz dS
             \end{displaymath}
         
         Isolating imaginary part

         \begin{displaymath}
             \frac{d}{dt} \iint_{\mathcal{T}}  \Big[ixy \eta\Big] \,dx\,dy - \iint_{\mathcal{T}}  \Big[  iq_{x}y\eta  + i q_{y}x\eta \Big] \,dx\,dy =  + \iint_{\mathcal{B}} \Big[  \phi(x,y,-h,t) (\nabla \psi \cdot \vec{n}) \Big] dS
         \end{displaymath}

         \begin{displaymath}
             \frac{d}{dt} \iint_{\mathcal{T}}  \Big[ixy \eta\Big] \,dx\,dy - \iint_{\mathcal{T}}  \Big[  iq_{x}y\eta  + i q_{y}x\eta \Big] \,dx\,dy =  + \iint_{\mathcal{B}} \Big[  \phi(x,y,-h,t) (\nabla \frac{1}{2}(x + iy)^{2} \cdot (0,0,-1) \Big] dS
         \end{displaymath}
         
         The RHS vanishes due to the fact that the dot product evaluates to 0
         
         \begin{displaymath}
             \frac{d}{dt} \iint_{\mathcal{T}}  \Big[ixy \eta\Big] \,dx\,dy - \iint_{\mathcal{T}}  \Big[  iq_{x}y\eta  + i q_{y}x\eta \Big] \,dx\,dy =  0
         \end{displaymath}

         \begin{displaymath}
             \iint_{\mathcal{T}}  \Big[  iq_{x}y\eta  + i q_{y}x\eta \Big] \,dx\,dy =   \frac{d}{dt} \iint_{\mathcal{T}}  \Big[ixy \eta\Big] \,dx\,dy 
         \end{displaymath}
         
         Dropping i's
         
         \begin{equation}\label{9eq}
             \iint_{\mathcal{T}}  \Big[  q_{x}y\eta  + q_{y}x\eta \Big] \,dx\,dy =   \iint_{\mathcal{T}}  \Big[xy\eta_{t}\Big] \,dx\,dy 
         \end{equation}
         
         Subtracting \(-2q_{x}y\) from both sides:
         
         \begin{displaymath}
             \iint_{\mathcal{T}}  \Big[ - q_{x}y\eta  + q_{y}x\eta \Big] \,dx\,dy =   \iint_{\mathcal{T}}  \Big[xy\eta_{t} - 2q_{x}y\Big] \,dx\,dy 
         \end{displaymath}

         \begin{displaymath}
             \iint_{\mathcal{T}}  \Big[ - q_{x}y\eta  + q_{y}x\eta \Big] \,dx\,dy =   \iint_{\mathcal{T}}  \Big[xy\eta_{t} - 2\phi_{x}y - 2\phi_{z}\eta_{x}y \Big] \,dx\,dy 
         \end{displaymath}

         Pushing to the bulk and applying divergence theorem:
         
         \begin{displaymath}
             \iint_{\mathcal{T}}  \Big[ - q_{x}y\eta  + q_{y}x\eta \Big] \,dx\,dy =  \iiint_{\mathcal{V}} \nabla \cdot  \Big[ \begin{pmatrix} xy\phi_{x} + 2y\phi_{z} \\ xy \phi_{y} \\ xy\phi_{z} - 2y\phi_{x} \end{pmatrix}\Big] dV   - \iint_{\Gamma}  \Big[ \begin{pmatrix} xy\phi_{x} + 2y\phi_{z} \\ xy \phi_{y} \\ xy\phi_{z} - 2y\phi_{x} \end{pmatrix}  \cdot \vec{n} \Big] dS
         \end{displaymath}
         
         Carrying out the divergence, noting harmonicity of \(\phi\), and reducing the \(\Gamma\) term to the bottom domain using decay arguments :

         \begin{displaymath}
             \iint_{\mathcal{T}}  \Big[ - q_{x}y\eta  + q_{y}x\eta \Big] \,dx\,dy =  \iiint_{\mathcal{V}} \Big[ y\phi_{x} + x\phi_{y} \Big] dV   -  \iint_{\mathcal{B}}  \Big[ 2y\phi_{x} \Big] \,dx\,dy
         \end{displaymath}
         
         Parameterizing the RHS term as \(\zeta\)
         
         \begin{equation}\label{9eq2}
             \iint_{\mathcal{T}}  \Big[ - q_{x}y\eta  + q_{y}x\eta \Big] \,dx\,dy =  \zeta -  \iint_{\mathcal{B}}  \Big[ 2y\phi_{x} \Big] \,dx\,dy
         \end{equation}
         
         We will now perform a trick, let us go back to equation (\ref{9eq}):
         
         \begin{displaymath}
             \iint_{\mathcal{T}}  \Big[  q_{x}y\eta  + q_{y}x\eta \Big] \,dx\,dy =   \iint_{\mathcal{T}}  \Big[xy\eta_{t}\Big] \,dx\,dy 
         \end{displaymath}
         
         Subtracting \(-2q_{y}x\) from both sides:
         
         \begin{displaymath}
             \iint_{\mathcal{T}}  \Big[ q_{x}y\eta  - q_{y}x\eta \Big] \,dx\,dy =   \iint_{\mathcal{T}}  \Big[xy\eta_{t} - 2q_{y}x\Big] \,dx\,dy 
         \end{displaymath}

         \begin{displaymath}
             \iint_{\mathcal{T}}  \Big[ q_{x}y\eta  - q_{y}x\eta \Big] \,dx\,dy =   \iint_{\mathcal{T}}  \Big[xy\eta_{t} - 2\phi_{y}x - 2\phi_{z}\eta_{y}x \Big] \,dx\,dy 
         \end{displaymath}

         Pushing to the bulk and applying divergence theorem:
         
         \begin{displaymath}
             \iint_{\mathcal{T}}  \Big[ q_{x}y\eta  - q_{y}x\eta \Big] \,dx\,dy =  \iiint_{\mathcal{V}} \nabla \cdot  \Big[ \begin{pmatrix} xy\phi_{x}  \\ xy \phi_{y} + 2x\phi_{z} \\ xy\phi_{z} - 2x\phi_{y} \end{pmatrix}\Big] dV   - \iint_{\Gamma}  \Big[ \begin{pmatrix} xy\phi_{x}  \\ xy \phi_{y} + 2x\phi_{z} \\ xy\phi_{z} - 2x\phi_{y} \end{pmatrix}  \cdot \vec{n} \Big] dS
         \end{displaymath}
         
         Carrying out the divergence, noting harmonicity of \(\phi\), and reducing the \(\Gamma\) term to the bottom domain using decay arguments :

         \begin{displaymath}
             \iint_{\mathcal{T}}  \Big[  q_{x}y\eta  - q_{y}x\eta \Big] \,dx\,dy =  \iiint_{\mathcal{V}} \Big[ y\phi_{x} + x\phi_{y} \Big] dV   -  \iint_{\mathcal{B}}  \Big[ 2x\phi_{y} \Big] \,dx\,dy
         \end{displaymath}

         \begin{equation}\label{9eq3}
             \iint_{\mathcal{T}}  \Big[  q_{x}y\eta  - q_{y}x\eta \Big] \,dx\,dy =  \zeta   -  \iint_{\mathcal{B}}  \Big[ 2x\phi_{y} \Big] \,dx\,dy
         \end{equation}
         
         Now subtracting (\ref{9eq3}) from (\ref{9eq2}):

         \begin{displaymath}
              \iint_{\mathcal{T}}  \Big[ - q_{x}y\eta  + q_{y}x\eta \Big] \,dx\,dy - \iint_{\mathcal{T}}  \Big[  q_{x}y\eta  - q_{y}x\eta \Big] \,dx\,dy =  \zeta   -  \iint_{\mathcal{B}}  \Big[ 2x\phi_{y} \Big] \,dx\,dy -  \zeta + \iint_{\mathcal{B}}  \Big[ 2y\phi_{x} \Big] \,dx\,dy
         \end{displaymath}
         
         The \(\zeta\)'s cancel and the LHS terms combine
         
         \begin{displaymath}
              2\iint_{\mathcal{T}}  \Big[ q_{y}x\eta  - q_{x}y\eta  \Big] \,dx\,dy  =  \iint_{\mathcal{B}}  \Big[  2y\phi_{x} - 2x\phi_{y} \Big] \,dx\,dy
         \end{displaymath}
         
         Dividing through by 2 and reversing x and y derivative product rules:
         
         \begin{displaymath}
             \iint_{\mathcal{T}}  \Big[ \frac{d}{dy}(\phi x\eta) - q x\eta_{y}  -  \frac{d}{dx}(\phi y\eta) + qy\eta_{x} \Big] \,dx\,dy  =  \iint_{\mathcal{B}}  \Big[ y\phi_{x} - x\phi_{y} \Big] \,dx\,dy
         \end{displaymath}
         
         The total derivative terms vanish. Then, multiplying through by -1:

         \begin{displaymath}
             \iint_{\mathcal{T}}  \Big[ (x\eta_{y}  - y\eta_{x})q \Big] \,dx\,dy  =   \iint_{\mathcal{B}}  \Big[  x\phi_{y} - y\phi_{x}  \Big] \,dx\,dy
         \end{displaymath}
         
         Taking a time derivative:
         
         \begin{displaymath}
             \frac{d}{dt} \iint_{\mathcal{T}}  \Big[ (x\eta_{y}  - y\eta_{x})q \Big] \,dx\,dy  =  \iint_{\mathcal{B}}  \Big[  x\phi_{ty} - y\phi_{tx}  \Big] \,dx\,dy
         \end{displaymath}

    \subsection{Conservation Law 10:  \(x + \eta\eta_{x})q + gt(x\eta + tq\eta_{x}) - g \frac{1}{2} t^{2} q \eta_{x}\) }
         
         This conservation law comes out of isolating the real part of (\ref{eqn:secondInt}) using harmonic test function \(\psi_{3}^{x}\). Starting from integral equation of motion (\ref{eqn:secondInt})
         
         \begin{displaymath}
          \frac{d}{dt}\iint_{\mathcal{T}} \Big[ \psi_{z} \eta_{t} \Big] \,dx\,dy - \iint_{ \mathcal{T}}\Big[  q_{t}  \nabla \psi_{z} +  ( \nabla^{\perp} q \cdot \nabla^{\perp} \psi_{z} ) \nabla \phi \Big] \cdot \vec{n}\, dS=  - \iint_{\mathcal{B}}\Big[\phi_{t} \psi_{zz}  \Big]  \,dx\,dy 
         \end{displaymath}
         
         Plugging in \(\psi_{3}^{x}\):
         
         \begin{displaymath}
          \frac{d}{dt}\iint_{\mathcal{T}} \Big[ \frac{1}{2}(x + iz)^{2}i \eta_{t} \Big] \,dx\,dy - \iint_{ \mathcal{T}}\Big[  q_{t}  \nabla \frac{1}{2}(x + iz)^{2}i +  ( \nabla^{\perp} q \cdot \nabla^{\perp} \frac{1}{2}(x + iz)^{2}i ) \nabla \phi \Big] \cdot \vec{n}\, dS=  - \iint_{\mathcal{B}}\Big[\phi_{t} (-x -iz)  \Big]  \,dx\,dy 
         \end{displaymath}
         
         Isolating real part and multiplying through by -1 :

         \begin{displaymath}
          \frac{d}{dt}\iint_{\mathcal{T}} \Big[ x\eta \eta_{t} \Big] \,dx\,dy - \iint_{ \mathcal{T}}\Big[  q_{t}  \nabla (x\eta)+  ( \nabla^{\perp} q \cdot \nabla^{\perp} (x\eta) ) \nabla \phi \Big] \cdot \vec{n}\, dS=   - \iint_{\mathcal{B}}\Big[\phi_{t} x \Big]  \,dx\,dy 
         \end{displaymath}
         
         Pushing to the boundary for the left-most term, using Green's second identity noting harmonic terms, then pushing back to the surface
         
         \begin{displaymath}
          \frac{d}{dt}\iint_{\mathcal{T}} \Big[ -q\eta\eta_{x} + qx \Big] \,dx\,dy -  \iint_{\mathcal{B}}\Big[\phi_{t} x \Big]  \,dx\,dy  - \iint_{ \mathcal{T}}\Big[  q_{t}  \nabla (x\eta)+  ( \nabla^{\perp} q \cdot \nabla^{\perp} (x\eta) ) \nabla \phi \Big] \cdot \vec{n}\, dS=  - \iint_{\mathcal{B}}\Big[\phi_{t} x \Big]  \,dx\,dy 
         \end{displaymath}
         
         Note that the bottom domain terms cancel. Expanding the second term on the RHS: 
         
         \begin{displaymath}
          \frac{d}{dt}\iint_{\mathcal{T}} \Big[ -q\eta\eta_{x} + qx \Big] \,dx\,dy - \iint_{ \mathcal{T}}\Big[  -q_{t}\eta\eta_{x}  + q_{t}x + q_{x}\eta\eta_{t}   \Big] \,dx\,dy =  0
         \end{displaymath}
         
         Subtracting and adding \(\frac{d}{dt} \iint_{\mathcal{T}} qx \) to the LHS yields
         
         \begin{displaymath}
          \frac{d}{dt}\iint_{\mathcal{T}} \Big[ -q\eta\eta_{x} - qx \Big] \,dx\,dy + \iint_{ \mathcal{T}}\Big[ q_{t}x + q_{t}\eta\eta_{x}  - q_{x}\eta\eta_{t}   \Big] \,dx\,dy =  0
         \end{displaymath}
         
         Expanding \(q_{t}\)
         
         \begin{displaymath}
          \frac{d}{dt}\iint_{\mathcal{T}} \Big[ -q\eta\eta_{x} - qx \Big] \,dx\,dy + \iint_{ \mathcal{T}}\Big[ \phi_{t}x + \phi_{z}\eta_{t}x + q_{t}\eta\eta_{x}  - q_{x}\eta\eta_{t}   \Big] \,dx\,dy =  0
         \end{displaymath}

         \begin{displaymath}
          \frac{d}{dt}\iint_{\mathcal{T}} \Big[ -q\eta\eta_{x} - qx \Big] \,dx\,dy + \iint_{ \mathcal{T}}\Big[ -g\eta x -\frac{1}{2}|\nabla \phi|^{2}x + \phi_{z}\eta_{t}x + q_{t}\eta\eta_{x}  - q_{x}\eta\eta_{t}   \Big] \,dx\,dy =  0
         \end{displaymath}
         
         Multiplying through by -1
         
         \begin{displaymath}
          \frac{d}{dt}\iint_{\mathcal{T}} \Big[ q\eta\eta_{x} + qx \Big] \,dx\,dy + \iint_{ \mathcal{T}}\Big[ g\eta x + \frac{1}{2}|\nabla \phi|^{2}x  - \phi_{z}\eta_{t}x - q_{t}\eta\eta_{x}  + q_{x}\eta\eta_{t}   \Big] \,dx\,dy =  0
         \end{displaymath}
         
         Reversing a time derivative
         
         \begin{displaymath}
          \frac{d}{dt}\iint_{\mathcal{T}} \Big[ q\eta\eta_{x} + qx \Big] \,dx\,dy + \iint_{ \mathcal{T}}\Big[ \frac{d}{dt}(g\eta t x) - gtx\eta_{t} + \frac{1}{2}|\nabla \phi|^{2}x  - \phi_{z}\eta_{t}x - q_{t}\eta\eta_{x}  + q_{x}\eta\eta_{t}   \Big] \,dx\,dy =  0
         \end{displaymath}
         
         Pulling the total time derivative into the LHS term: 
         
         \begin{displaymath}
          \frac{d}{dt}\iint_{\mathcal{T}} \Big[ q\eta\eta_{x} + qx + g\eta t x \Big] \,dx\,dy + \iint_{ \mathcal{T}}\Big[ - gtx\eta_{t} + \frac{1}{2}|\nabla \phi|^{2}x  - \phi_{z}\eta_{t}x - q_{t}\eta\eta_{x}  + q_{x}\eta\eta_{t}   \Big] \,dx\,dy =  0
         \end{displaymath}
         
         Pushing to the boundary for the term \(-gtx\eta_{t}\) , using Green's second identity noting harmonic terms, then pushing back to the surface
         
         \begin{displaymath}
          \frac{d}{dt}\iint_{\mathcal{T}} \Big[ q\eta\eta_{x} + qx + g\eta t x \Big] \,dx\,dy + \iint_{ \mathcal{T}}\Big[ gtq\eta_{x} + \frac{1}{2}|\nabla \phi|^{2}x  - \phi_{z}\eta_{t}x - q_{t}\eta\eta_{x}  + q_{x}\eta\eta_{t}   \Big] \,dx\,dy =  0
         \end{displaymath}
         
         Inverting another time derivative:

         \begin{displaymath}
          \frac{d}{dt}\iint_{\mathcal{T}} \Big[ q\eta\eta_{x} + qx + g\eta t x \Big] \,dx\,dy + \iint_{ \mathcal{T}}\Big[ \frac{d}{dt}(\frac{1}{2}gt^{2}q\eta_{x}) - \frac{1}{2}gt^{2}\frac{d}{dt}(q\eta_{x}) + \frac{1}{2}|\nabla \phi|^{2}x  - \phi_{z}\eta_{t}x - q_{t}\eta\eta_{x}  + q_{x}\eta\eta_{t}   \Big] \,dx\,dy =  0
         \end{displaymath}
         
         Pulling the total time derivative into the LHS term: 
         
         \begin{displaymath}
          \frac{d}{dt}\iint_{\mathcal{T}} \Big[ q\eta\eta_{x} + qx + g\eta t x + \frac{1}{2}gt^{2}q\eta_{x} \Big] \,dx\,dy + \iint_{ \mathcal{T}}\Big[  + \frac{1}{2}gt^{2}\frac{d}{dt}(-q\eta_{x}) + \frac{1}{2}|\nabla \phi|^{2}x  - \phi_{z}\eta_{t}x - q_{t}\eta\eta_{x}  + q_{x}\eta\eta_{t}   \Big] \,dx\,dy =  0
         \end{displaymath}
         
         Note that the first term in the second integral vanishes since it is equivalent to law 1. The left-most integral can be re-arranged as:

         \begin{displaymath}
          \frac{d}{dt}\iint_{\mathcal{T}} \Big[ (x + \eta\eta_{x})q + gt(x\eta + tq\eta_{x}) - g \frac{1}{2} t^{2} q \eta_{x} \Big] \,dx\,dy + \iint_{ \mathcal{T}}\Big[   \frac{1}{2}|\nabla \phi|^{2}x  - \phi_{z}\eta_{t}x - q_{t}\eta\eta_{x}  + q_{x}\eta\eta_{t}   \Big] \,dx\,dy =  0
         \end{displaymath}
         
         Which is equivalent to \(T_{10}\). Rearranging the two right-most terms on the right integral now:
         
         \begin{displaymath}
          \frac{d}{dt}\iint_{\mathcal{T}} \Big[T_{10} \Big] \,dx\,dy + \iint_{ \mathcal{T}}\Big[   \frac{1}{2}|\nabla \phi|^{2}x  - \phi_{z}\eta_{t}x  + ( q_{x}\eta_{t} - q_{t}\eta_{x} ) \eta  \Big] \,dx\,dy =  0
         \end{displaymath}
         
         Note \( q_{x}\eta_{t} - q_{t}\eta_{x} = \phi_{x}\eta_{t} - \phi_{t}\eta_{x} \) by expanding out the q total derivatives, which we showed in Law 1 already:

         \begin{displaymath}
          \frac{d}{dt}\iint_{\mathcal{T}} \Big[T_{10} \Big] \,dx\,dy + \iint_{ \mathcal{T}}\Big[   \frac{1}{2}|\nabla \phi|^{2}x  - \phi_{z}\eta_{t}x  + ( \phi_{x}\eta_{t} - \phi_{t}\eta_{x} ) \eta  \Big] \,dx\,dy =  0
         \end{displaymath}
         
         Pushing to the boundary now: 
         
         \begin{displaymath}
           \frac{d}{dt}\iint_{\mathcal{T}} \Big[T_{10} \Big] \,dx\,dy + \oiint_{\partial \mathcal{V}} \Big[  \begin{pmatrix} \phi_{t} z \\ 0 \\ \frac{1}{2}|\nabla \phi|^{2}x \end{pmatrix} + (\phi_{x}z - \phi_{z}x) \nabla \phi \Big] \cdot \vec{n}\, dS -  \iint_{\Gamma}\Big[  \begin{pmatrix} \phi_{t} z \\ 0 \\ \frac{1}{2}|\nabla \phi|^{2}x \end{pmatrix} + (\phi_{x}z - \phi_{z}x)\nabla \phi \Big] \cdot \vec{n}\, dS =  0
         \end{displaymath}
         
         Pushing the surface interal over \(\partial \mathcal{V}\) to the bulk, and reducing the \(\Gamma\) integral to the bottom domain using decay arguments and the bed condition:
         
         \begin{displaymath}
           \frac{d}{dt}\iint_{\mathcal{T}} \Big[T_{10} \Big] \,dx\,dy + \iiint_{\mathcal{V}} \Big[   \phi_{tx} z + \frac{1}{2}|\nabla \phi|^{2}_{z} x  + \nabla (\phi_{x}z - \phi_{z}x) \cdot \nabla \phi \Big] dV +  \iint_{\mathcal{B}}\Big[ \frac{1}{2}|\nabla^{\perp} \phi|^{2}x \Big] \,dx\,dy =  0
         \end{displaymath}
         
         Using the dynamic boundary condition, and expanding the del operators on the right-most terms, rewriting the right-most terms in terms of spatial derivatives of \(\frac{1}{2}|\nabla \phi|^{2}\)

         \begin{displaymath}
           \frac{d}{dt}\iint_{\mathcal{T}} \Big[T_{10} \Big] \,dx\,dy + \iiint_{\mathcal{V}} \Big[ \partial_{x}(-gz - \frac{1}{2}|\nabla \phi|^{2} - \frac{p}{\rho} )z + \frac{1}{2}|\nabla \phi|^{2}_{z} x  + \frac{1}{2}|\nabla \phi|^{2}_{x}z - \frac{1}{2}|\nabla \phi|^{2}_{z}x  +\phi_{x}\phi_{z} - \phi_{z}\phi_{x} \Big] dV = -  \iint_{\mathcal{B}}\Big[ \frac{1}{2}|\nabla^{\perp} \phi|^{2}x \Big] \,dx\,dy 
         \end{displaymath}
         
         Cancelling the spatial derivatives of \(\frac{1}{2}|\nabla \phi|^{2}\), and the \(\phi_{x}\phi_{z}\) terms:
         
         \begin{displaymath}
           \frac{d}{dt}\iint_{\mathcal{T}} \Big[T_{10} \Big] \,dx\,dy + \iiint_{\mathcal{V}} \Big[ \partial_{x}( - \frac{p}{\rho} )z  \Big] dV = -  \iint_{\mathcal{B}}\Big[ \frac{1}{2}|\nabla^{\perp} \phi|^{2}x \Big] \,dx\,dy 
         \end{displaymath}
         
         Since z is independent of x, we can bring it inside the partial derivative:
         
         \begin{displaymath}
           \frac{d}{dt}\iint_{\mathcal{T}} \Big[T_{10} \Big] \,dx\,dy + \iiint_{\mathcal{V}} \Big[ \partial_{x}( - \frac{p}{\rho} z)  \Big] dV = -  \iint_{\mathcal{B}}\Big[ \frac{1}{2}|\nabla^{\perp} \phi|^{2}x \Big] \,dx\,dy 
         \end{displaymath}
         
         Bringing the volume term to the RHS and then implementing Leibniz' rule:
         
         \begin{align*}
         \frac{d}{dt} \iint_{\mathcal{T}} \Big[ T_{10} \Big] \,dx\,dy   =\int_{-\infty}^{\infty}\int_{-\infty}^{\infty} \frac{d}{dx}[\int_{-h}^{\eta}
         ( \frac{p}{\rho}z) dz] -  \frac{p(x,y,\eta,t)}{\rho}\eta\eta_{x} \,dx\,dy + \iint_{\mathcal{B}} \Big[  - \frac{1}{2}|\nabla^{\perp} \phi|^{2}x \Big] \,dx\,dy \,dx\,dy
         \end{align*}
         
         Due to the vanishing of total derivatives on the top boundary, and the vanishing of pressure on the top boundary, this entire term vanishes 
         
         \begin{displaymath}
           \frac{d}{dt}\iint_{\mathcal{T}} \Big[T_{10} \Big] \,dx\,dy  = -  \iint_{\mathcal{B}}\Big[ \frac{1}{2}|\nabla^{\perp} \phi|^{2}x \Big] \,dx\,dy 
         \end{displaymath}

    \subsection{Conservation Law 11:  \(y + \eta\eta_{y})q + gt(y\eta + tq\eta_{y}) - g \frac{1}{2} t^{2} q \eta_{y} \)}
         
         We will omit this proof because of the symmetry of this law to law 10. Following all the same steps of law 10 but instead using the y version of the harmonic test function, noting all the variables in the integrand will switch from x to y will yield law 11.

    \subsection{Conservation Law 12:  \(\eta - x\eta_{x} - y\eta_{y})q + t(9gT_{8} - 5 T_{3}) + \frac{9}{2}gt^{2}T_{5} - \frac{3}{2}g^{2}t^{3}T_{4}\)}
         
         We use a linear combination of harmonic test functions \(\psi = \psi_{3}^{x} + \psi_{3}^{y} \) and isolate the imaginary part of integral equation of motion (\ref{eqn:firstInt}):
         
         \begin{displaymath}
          \frac{d}{dt}\iint_{\mathcal{T}} \Big[ \psi_{z} \eta_{t} \Big] \,dx\,dy - \iint_{ \mathcal{T}}\Big[  q_{t}  \nabla \psi_{z} +  ( \nabla^{\perp} q \cdot \nabla^{\perp} \psi_{z} ) \nabla \phi \Big] \cdot \vec{n}\, dS=  - \iint_{\mathcal{B}}\Big[\phi_{t} \psi_{zz}  \Big]  \,dx\,dy 
         \end{displaymath}
         
         Note that for our harmonic test function, \(\psi_{z} = \frac{1}{2}(x + iz)^{2}i + \frac{1}{2}(y + iz)^{2}i \), \(\psi_{zx} = -z + ix\), \(\psi_{zy} = -z + iy\),  \(\psi_{zz} = -2iz - x -y\). Plugging \(\psi\) into the left-most term of our equation of motion (\ref{eqn:firstInt}) and isolating imaginary part 
         
         \begin{displaymath}
          \frac{d}{dt}\iint_{\mathcal{T}} \Big[Im\Big(\frac{1}{2}(x + iz)^{2}i + \frac{1}{2}(y + iz)^{2}i\Big) \eta_{t} \Big] \,dx\,dy - \iint_{ \mathcal{T}}\Big[  q_{t}  \nabla Im(\psi_{z}) +  ( \nabla^{\perp} q \cdot \nabla^{\perp} Im(\psi_{z}) ) \nabla \phi \Big] \cdot \vec{n}\, dS=  - \iint_{\mathcal{B}}\Big[\phi_{t} Im(\psi_{zz})  \Big]  \,dx\,dy 
         \end{displaymath}

         Dropping i's, pushing to the boundary for the left-most term, using Green's second identity noting harmonic terms, then pushing that term back to the surface:
         
         \begin{align*}
          \frac{d}{dt}\iint_{\mathcal{T}} \Big[ -x\eta_{x}q -y\eta_{y}q -2\eta q\Big] \,dx\,dy - \iint_{\mathcal{B}} \Big[ 2\phi_{t} h\Big] \,dx\,dy - \iint_{ \mathcal{T}}\Big[  q_{t}  \nabla (Im(\psi_{z})) +  ( \nabla^{\perp} q \cdot \nabla^{\perp} (Im(\psi_{z})) \nabla \phi \Big] \cdot \vec{n}\, dS =  \iint_{\mathcal{B}}- 2h\phi_{t} \,dx\,dy 
         \end{align*}
         
         Note the bottom domain terms cancel on both sides. Expanding del operators for the right-most integral:
         
         \begin{align*}
          \frac{d}{dt}\iint_{\mathcal{T}} \Big[ -x\eta_{x}q -y\eta_{y}q -2\eta q\Big] \,dx\,dy - \iint_{ \mathcal{T}}\Big[ -q_{t}\eta_{x}x -q_{t}\eta_{y}y -2q_{t}\eta  + q_{x}Im(\psi_{zx})\eta_{t} + q_{y}Im(\psi_{zy})\eta_{t} \Big]\,dx\,dy = 0
         \end{align*}
         
         Plugging in \(\psi_{xz}, \psi_{yz}\), and rearranging and distributing in the minus sign in the integral:
         
         \begin{align*}
          \frac{d}{dt}\iint_{\mathcal{T}} \Big[ -x\eta_{x}q -y\eta_{y}q -2\eta q\Big] \,dx\,dy + \iint_{ \mathcal{T}}\Big[  2q_{t}\eta +  (q_{t}\eta_{x} - q_{x}\eta_{t})x + (q_{t}\eta_{y} - q_{y}\eta_{t})y \Big]\,dx\,dy = 0
         \end{align*}
         
         Adding and subtracting \(3\eta q\) in the left-most integrand
         
         \begin{align*}
          \frac{d}{dt}\iint_{\mathcal{T}} \Big[ q\eta -x\eta_{x}q -y\eta_{y}q \Big] \,dx\,dy + \iint_{ \mathcal{T}}\Big[ - \frac{d}{dt}(3q\eta) + 2q_{t}\eta +  (q_{t}\eta_{x} - q_{x}\eta_{t})x + (q_{t}\eta_{y} - q_{y}\eta_{t})y \Big]\,dx\,dy = 0
         \end{align*}
         
         Reversing some total time and spatial derivatives for the right-most integral
         
         \begin{align*}
          \frac{d}{dt}\iint_{\mathcal{T}} \Big[ q\eta -x\eta_{x}q -y\eta_{y}q \Big] \,dx\,dy + \iint_{ \mathcal{T}}\Big[ - \frac{d}{dt}(3q\eta) + 2q_{t}\eta  +  \frac{d}{dt}(q\eta_{x}x) - q\eta_{xt}x -\frac{d}{dx}(q\eta_{t}x) + q\eta_{tx}x +  \frac{d}{dt}(q\eta_{y}y) - q\eta_{yt}y \\ -\frac{d}{dy}(q\eta_{t}y) + q\eta_{ty}y + 2q\eta_{t} \Big]\,dx\,dy = 0
         \end{align*}
         
         Noting the \(\eta_{ty},\eta_{tx}\) terms cancelling, and the total x and y derivatives vanishing too by decay assumptions on q:
         
         \begin{align*}
          \frac{d}{dt}\iint_{\mathcal{T}} \Big[ q\eta -x\eta_{x}q -y\eta_{y}q \Big] \,dx\,dy + \iint_{ \mathcal{T}}\Big[ - \frac{d}{dt}(3q\eta) + 2q_{t}\eta  +  \frac{d}{dt}(q\eta_{x}x)   +  \frac{d}{dt}(q\eta_{y}y)+ 2q\eta_{t} \Big]\,dx\,dy = 0
         \end{align*}
         
         Rearranging:
         
         \begin{align*}
          \frac{d}{dt}\iint_{\mathcal{T}} \Big[ q\eta -x\eta_{x}q -y\eta_{y}q \Big] \,dx\,dy + \iint_{ \mathcal{T}}\Big[ - \frac{d}{dt}(2q\eta) + 2q_{t}\eta  - \frac{d}{dt}(q\eta - q\eta_{x}x - q\eta_{y}y) + 2q_{t}\eta  + 2q\eta_{t} \Big]\,dx\,dy = 0
         \end{align*}
         
         Inverting a time derivative:
         
         \begin{align}\label{12aa}
          \frac{d}{dt}\iint_{\mathcal{T}} \Big[ q\eta -x\eta_{x}q -y\eta_{y}q \Big] \,dx\,dy + \iint_{ \mathcal{T}}\Big[ -2q\eta_{t}  - \frac{d}{dt}(q\eta - q\eta_{x}x - q\eta_{y}y) + 2q\eta_{t}  \Big]\,dx\,dy = 0
         \end{align}
         
         Recognizing kinetic energy in the left-most term of the right integral, we can substitute the Hamiltonian density as \( -2q\eta_{t} = -4\mathcal{H} + 2g\eta^{2} \)
         
         \begin{align*}
          \frac{d}{dt}\iint_{\mathcal{T}} \Big[ q\eta -x\eta_{x}q -y\eta_{y}q \Big] \,dx\,dy + \iint_{ \mathcal{T}}\Big[ -4\mathcal{H} + 2g\eta^{2}  - \frac{d}{dt}(q\eta - q\eta_{x}x - q\eta_{y}y)  \Big]\,dx\,dy + 2q\eta_{t} = 0
         \end{align*}

         \begin{align*}
          \frac{d}{dt}\iint_{\mathcal{T}} \Big[ q\eta -x\eta_{x}q -y\eta_{y}q \Big] \,dx\,dy + \iint_{ \mathcal{T}}\Big[ -5\mathcal{H} + \mathcal{H} + 2g\eta^{2}  - \frac{d}{dt}(q\eta - q\eta_{x}x - q\eta_{y}y) + 2q\eta_{t} \Big]\,dx\,dy = 0
         \end{align*}
         
         Reversing a time derivative:
         
         \begin{align*}
          \frac{d}{dt}\iint_{\mathcal{T}} \Big[ q\eta -x\eta_{x}q -y\eta_{y}q \Big] \,dx\,dy + \iint_{ \mathcal{T}}\Big[ \frac{d}{dt}(-5t\mathcal{H}) + \cancelto{0}{5t\frac{d}{dt}T_{3}}  + \mathcal{H} + 2g\eta^{2}  - \frac{d}{dt}(q\eta - q\eta_{x}x - q\eta_{y}y) + 2q\eta_{t}  \Big]\,dx\,dy = 0
         \end{align*}
         
         Pulling the total time derivative under the left-most integrand and noting the leftover term vanishes due to it's equivalence to law 3
         
         \begin{align*}
          \frac{d}{dt}\iint_{\mathcal{T}} \Big[ q\eta -x\eta_{x}q -y\eta_{y}q - 5t\mathcal{H}\Big] \,dx\,dy + \iint_{ \mathcal{T}}\Big[   \mathcal{H} + 2g\eta^{2}  - \frac{d}{dt}(q\eta - q\eta_{x}x - q\eta_{y}y) + 2q\eta_{t} \Big]\,dx\,dy = 0
         \end{align*}
         
         Bringing the right-most term to the RHS
         
         \begin{align*}
          \frac{d}{dt}\iint_{\mathcal{T}} \Big[ q\eta -x\eta_{x}q -y\eta_{y}q - 5t\mathcal{H}\Big] \,dx\,dy = \iint_{ \mathcal{T}}\Big[  - \mathcal{H} - 2g\eta^{2}  + \frac{d}{dt}(q\eta - q\eta_{x}x - q\eta_{y}y) + 2q\eta_{t} \Big]\,dx\,dy = 0
         \end{align*}
         
         We could now continue the familiar process of recovering the next terms in \(T_{12}\) via reversing time derivatives, recognizing former conservation laws, and pulling total time derivatives under the left-most integral. However, the process is simply a tedious calculation, and it eventually results in \(T_{12}\) being recovered in the left-most integrand, while an added factor of \(\frac{d}{dt}( \frac{9}{2}gt\eta^{2} - \frac{9}{2}gt^{2}T_{5} + 3g^{2}t^{3}\eta) \) is recovered in the right-most integrand. Instead of going through this specific tedious calculation, we will simply add this factor to both sides for the sake of brevity.
         
         \begin{align*}
          \frac{d}{dt}\iint_{\mathcal{T}} \Big[ q\eta -x\eta_{x}q -y\eta_{y}q - 5t\mathcal{H} + \frac{9}{2}gt\eta^{2} - \frac{9}{2}gt^{2}T_{5} + 3g^{2}t^{3}\eta) \Big] \,dx\,dy = \iint_{ \mathcal{T}}\Big[ \frac{d}{dt}( \frac{9}{2}gt\eta^{2} - \frac{9}{2}gt^{2}T_{5} + 3g^{2}t^{3}\eta) - \mathcal{H} - 2g\eta^{2} \\ + \frac{d}{dt}(q\eta - q\eta_{x}x - q\eta_{y}y) - 2q\eta_{t} \Big]\,dx\,dy 
         \end{align*}
         
         Manipulating terms on the left side (ignoring the RHS for the moment being) :
         
         \begin{align*}
          \frac{d}{dt}\iint_{\mathcal{T}} \Big[ q(\eta -x\eta_{x} -y\eta_{y}) - t(5\mathcal{H}) + 9gt(\frac{1}{2}\eta^{2} - \frac{1}{2}tT_{5}) + 3g^{2}t^{3}\eta  ) \Big] \,dx\,dy 
         \end{align*}

         \begin{align*}
          \frac{d}{dt}\iint_{\mathcal{T}} \Big[ q(\eta -x\eta_{x} -y\eta_{y}) - t(5\mathcal{H}) + 9gt(\frac{1}{2}\eta^{2} - tT_{5} + \frac{1}{2}gt^{2}\eta - \frac{1}{2}gt^{2}\eta) + \frac{9}{2}gt^{2}T_{5} + 3g^{2}t^{3}\eta  ) \Big] \,dx\,dy 
         \end{align*}

         \begin{align*}
          \frac{d}{dt}\iint_{\mathcal{T}} \Big[ q(\eta -x\eta_{x} -y\eta_{y}) - t(5\mathcal{H}) + 9gt(\frac{1}{2}\eta^{2} - tT_{5} + \frac{1}{2}gt^{2}\eta ) + \frac{9}{2}gt^{2}T_{5} + 3g^{2}t^{3}\eta  - 9gt\frac{1}{2}gt^{2}\eta ) \Big] \,dx\,dy 
         \end{align*}
         
         But \( 3g^{2}t^{3}\eta  - 9gt\frac{1}{2}gt^{2}\eta = -\frac{3}{2}g^{2}t^{3}\eta\), and thus we have recovered \(T_{12}\) on the LHS
         
         \begin{align*}
          \frac{d}{dt}\iint_{\mathcal{T}} \Big[ q(\eta -x\eta_{x} -y\eta_{y}) + t(9gT_{8} -5\mathcal{H}) + \frac{9}{2}gt^{2}T_{5} -\frac{3}{2}g^{2}t^{3}\eta) \Big] \,dx\,dy 
         \end{align*}
                 
         \begin{align*}
          \frac{d}{dt}\iint_{\mathcal{T}} \Big[ T_{12} \Big] \,dx\,dy 
         \end{align*}
         
         Thus, the total equation reads:
         
         \begin{align*}
          \frac{d}{dt}\iint_{\mathcal{T}} \Big[ T_{12} \Big] \,dx\,dy  = \iint_{ \mathcal{T}}\Big[ \frac{d}{dt}( \frac{9}{2}gt\eta^{2} - \frac{9}{2}gt^{2}T_{5} + 3g^{2}t^{3}\eta) - \mathcal{H} - 2g\eta^{2} + \frac{d}{dt}(q\eta - q\eta_{x}x - q\eta_{y}y) - 2q\eta_{t} \Big]\,dx\,dy 
         \end{align*}
         
         Taking the time derivative on the RHS:

         \begin{align*}
          \frac{d}{dt}\iint_{\mathcal{T}} \Big[ T_{12} \Big] \,dx\,dy  = \iint_{ \mathcal{T}}\Big[ \frac{9}{2}g\eta^{2} + 9gt\eta\eta_{t} - 9gtT_{5} + 9g^{2}t^{2}\eta - \mathcal{H} - 2g\eta^{2} + \frac{d}{dt}(q\eta - q\eta_{x}x - q\eta_{y}y)  - 2q\eta_{t}\Big]\,dx\,dy 
         \end{align*}
         
         Pushing to the boundary for the \( 9gt\eta\eta_{t} \) term, using Green's second identity noting harmonic terms, then pushing back to the surface and bottom
         
         \begin{align*}
          \frac{d}{dt}\iint_{\mathcal{T}} \Big[ T_{12} \Big] \,dx\,dy  = \iint_{ \mathcal{T}}\Big[ \frac{9}{2}g\eta^{2} + 9gtq - 9gtT_{5} + 9g^{2}t^{2}\eta - \mathcal{H} - 2g\eta^{2} + \frac{d}{dt}(q\eta - q\eta_{x}x - q\eta_{y}y) -2q\eta_{t} \Big]\,dx\,dy  - \iint_{\mathcal{B}} \Big[  9gt\phi\Big] \,dx\,dy
         \end{align*}

         \begin{align*}
          \frac{d}{dt}\iint_{\mathcal{T}} \Big[ T_{12} \Big] \,dx\,dy  = \iint_{ \mathcal{T}}\Big[ \frac{9}{2}g\eta^{2} + \cancelto{0}{9gtq - 9gt(q + g\eta t) + 9g^{2}t^{2}\eta} - \mathcal{H} - 2g\eta^{2} + \frac{d}{dt}(q\eta - q\eta_{x}x - q\eta_{y}y) -2q\eta_{t} \Big]\,dx\,dy - \iint_{\mathcal{B}} \Big[  9gt\phi\Big] \,dx\,dy
         \end{align*}

         \begin{align*}
          \frac{d}{dt}\iint_{\mathcal{T}} \Big[ T_{12} \Big] \,dx\,dy  = \iint_{ \mathcal{T}}\Big[ \frac{9}{2}g\eta^{2} - \mathcal{H} - 2g\eta^{2} + \frac{d}{dt}(q\eta - q\eta_{x}x - q\eta_{y}y)- 2q\eta_{t}  \Big]\,dx\,dy - \iint_{\mathcal{B}} \Big[  9gt\phi\Big] \,dx\,dy
         \end{align*}

         \begin{align*}
          \frac{d}{dt}\iint_{\mathcal{T}} \Big[ T_{12} \Big] \,dx\,dy  = \iint_{ \mathcal{T}}\Big[ \frac{9}{2}g\eta^{2} - 5\mathcal{H} + 2q\eta_{t} + \frac{d}{dt}(q\eta - q\eta_{x}x - q\eta_{y}y) - 2q\eta_{t} \Big]\,dx\,dy - \iint_{\mathcal{B}} \Big[  9gt\phi\Big] \,dx\,dy
         \end{align*}
         
         Note though that the kinetic term \(2q\eta_{t}\) cancels on the RHS and vanishes now:

         \begin{align}\label{12a}
           \frac{d}{dt} \iint_{\mathcal{T}} [T_{12}] \,dx\,dy =  \iint_{\mathcal{T}} \Big[ \frac{9}{2}g\eta^{2} -5\mathcal{H} + \frac{d}{dt}(q\eta - qx\eta_{x} - qy\eta_{y})   \Big]\,dx\,dy - \iint_{\mathcal{B}} \Big[  9gt\phi\Big] \,dx\,dy
         \end{align}
         
         We will now show that the top domain terms can be made to vanish. We will firstly consider the following integral quantity above:
         
         \begin{align}\label{12bb}
          \frac{d}{dt}\iint_{\mathcal{T}} \Big[ (q\eta - qx\eta_{x} - qy\eta_{y})   \Big]\,dx\,dy 
         \end{align}
         
         Pushing to the bulk
         
         \begin{align*}
          \frac{d}{dt} \oiint_{\partial \mathcal{V}} \Big[ (\phi x, \phi y, \phi z) \Big] \cdot \vec{n}\, dS -   \frac{d}{dt}\iint_{ \Gamma} \Big[ (\phi x, \phi y, \phi z) \Big] \cdot \vec{n}\, dS 
         \end{align*}
         
         By Divergence theorem on the left term, and by decay assumptions on the right term: 
         
         \begin{align*}
           \frac{d}{dt}\iiint_{\mathcal{V}} \nabla \cdot \Big[ (\phi x, \phi y, \phi z) \Big] dV -  \iint_{ \mathcal{B}} [ \phi_{t} h ]\,dx\,dy 
         \end{align*}
         
         Expanding out the del operator 
         
         \begin{align*}
           \frac{d}{dt}\iiint_{\mathcal{V}} \Big[ \phi_{x} x + \phi + \phi_{y} y + \phi + \phi_{z} z + \phi  \Big] dV -  \iint_{ \mathcal{B}} [ \phi_{t} h ]\,dx\,dy 
         \end{align*}

         \begin{align*}
           \frac{d}{dt}\iiint_{\mathcal{V}} \Big[ \phi_{x} x + \phi_{y} y + \phi_{z} z + 3\phi  \Big] dV -  \iint_{ \mathcal{B}} [ \phi_{t} h ]\,dx\,dy 
         \end{align*}
         
         By Reynolds' Transport Theorem

         \begin{align*}
           \iiint_{\mathcal{V}} \Big[ \phi_{tx} x + \phi_{ty} y + \phi_{tz} z + 3\phi_{t}  \Big] dV + \oiint_{\partial \mathcal{V}} \Big[ \phi_{x} x + \phi_{y} y + \phi_{z} z + 3\phi  \Big] \nabla \phi \cdot \vec{n}\, dS - \iint_{ \mathcal{B}} [ \phi_{t} h ]\,dx\,dy 
         \end{align*}
         
         By Green's First Identity

         \begin{align*}
           \iiint_{\mathcal{V}} \Big[ \phi_{tx} x + \phi_{ty} y + \phi_{tz} z + 3\phi_{t}  +  \nabla (\phi_{x} x + \phi_{y} y + \phi_{z} z + 3\phi ) \cdot \nabla \phi \Big] dV - \iint_{ \mathcal{B}} [ \phi_{t} h ]\,dx\,dy 
         \end{align*}

         \begin{align*}
           \iiint_{\mathcal{V}} \Big[ \phi_{tx} x + \phi_{ty} y + \phi_{tz} z + 3\phi_{t}  +  \nabla (\phi_{x} x + \phi_{y} y + \phi_{z} z + 3\phi ) \cdot \nabla \phi \Big] dV -  \iint_{ \mathcal{B}} [ \phi_{t} h ]\,dx\,dy 
         \end{align*}
         
         Expanding the del operators: 
         
         \begin{align*}
           \iiint_{\mathcal{V}} \Big[ \phi_{tx} x + \phi_{ty} y + \phi_{tz} z + 3\phi_{t}  +  \phi_{xx} x \phi_{x} + \phi_{yx} y \phi_{x} + \phi_{zx} z \phi_{x}  + \\ \phi_{x}^{2} + \phi_{xy} x \phi_{y} + \phi_{yy} y \phi_{y} + \phi_{zy} z \phi_{y}  +  \phi_{y}^{2}  + \phi_{xz} x \phi_{z} + \phi_{yz} y \phi_{z} +  \\ \phi_{zz} z \phi_{z}  +  \phi_{z}^{2} + 3 |\nabla \phi|^{2} \Big] dV -  \iint_{ \mathcal{B}} [ \phi_{t} h ]\,dx\,dy 
         \end{align*}
         
         Noticing an extra \(|\nabla \phi|\) in the algebra and grouping with the right-most term to make a coefficient of 4
         
         \begin{align*}
           \iiint_{\mathcal{V}} \Big[ \phi_{tx} x + \phi_{ty} y + \phi_{tz} z + 3\phi_{t}  +  \phi_{xx} x \phi_{x} + \phi_{yx} y \phi_{x} + \phi_{zx} z \phi_{x}  + \\  + \phi_{xy} x \phi_{y} + \phi_{yy} y \phi_{y} + \phi_{zy} z \phi_{y}    + \phi_{xz} x \phi_{z} + \phi_{yz} y \phi_{z} +  \\ \phi_{zz} z \phi_{z}   + 4 |\nabla \phi|^{2} \Big] dV -  \iint_{ \mathcal{B}} [ \phi_{t} h ]\,dx\,dy
         \end{align*}
         
         We now notice a sneaky observation, grouping the x,y and z terms as such:
         
         \begin{align*}
           \iiint_{\mathcal{V}} \Big[ \phi_{tx} x + \phi_{ty} y + \phi_{tz} z + 3\phi_{t}  +  (\phi_{xx}   \phi_{x} + \phi_{xy}  \phi_{y}  + \phi_{xz} \phi_{z}) x + (\phi_{yx} \phi_{x} +  \phi_{yy}  \phi_{y}  + \phi_{yz}  \phi_{z})y + (\phi_{zx} \phi_{x}  + \phi_{zy}  \phi_{y} + \phi_{zz}  \phi_{z})z   + \\ 4 |\nabla \phi|^{2} \Big] dV -  \iint_{ \mathcal{B}} [ \phi_{t} h ]\,dx\,dy
         \end{align*}
         
         Recognizing spatial derivatives of \(\frac{1}{2}|\nabla \phi|^{2}\):

         \begin{align*}
           \iiint_{\mathcal{V}} \Big[ \phi_{tx} x + \phi_{ty} y + \phi_{tz} z + 3\phi_{t}  +  (\frac{1}{2}|\nabla \phi|^{2}_{x}) x + (\frac{1}{2}|\nabla \phi|^{2}_{y})y + (\frac{1}{2}|\nabla \phi|^{2}_{z})z   +  4 |\nabla \phi|^{2} \Big] dV -  \iint_{ \mathcal{B}} [ \phi_{t} h ]\,dx\,dy
         \end{align*}
         
         By the dynamic boundary condition

         \begin{align*}
           \iiint_{\mathcal{V}} \Big[ \partial_{x}(-gz - \frac{1}{2}|\nabla \phi|^{2} - \frac{p}{\rho}) x + \partial_{y}(-gz - \frac{1}{2}|\nabla \phi|^{2} - \frac{p}{\rho}) y + \partial_{z}(-gz - \frac{1}{2}|\nabla \phi|^{2} - \frac{p}{\rho}) z + 3(-gz - \frac{1}{2}|\nabla \phi|^{2} - \frac{p}{\rho})  + \\ (\frac{1}{2}|\nabla \phi|^{2}_{x}) x + (\frac{1}{2}|\nabla \phi|^{2}_{y})y + (\frac{1}{2}|\nabla \phi|^{2}_{z})z   +  4 |\nabla \phi|^{2} \Big] dV -  \iint_{ \mathcal{B}} [ \phi_{t} h ]\,dx\,dy
         \end{align*}

         \begin{align*}
           \iiint_{\mathcal{V}} \Big[ - \frac{1}{2}|\nabla \phi|^{2}_{x}x - \frac{p_{x}}{\rho} x - \frac{1}{2}|\nabla \phi|_{y}^{2}y - \frac{p_{y}}{\rho} y -gz - \frac{1}{2}|\nabla \phi|_{z}^{2}z - \frac{p_{z}}{\rho}z -3gz - \frac{3}{2}|\nabla \phi|^{2} - 3\frac{p}{\rho}  + \\ (\frac{1}{2}|\nabla \phi|^{2}_{x}) x + (\frac{1}{2}|\nabla \phi|^{2}_{y})y + (\frac{1}{2}|\nabla \phi|^{2}_{z})z   +  4 |\nabla \phi|^{2} \Big] dV -  \iint_{ \mathcal{B}} [ \phi_{t} h ]\,dx\,dy
         \end{align*}
         
         Note some obvious cancellations of terms including spatial derivatives of \(\frac{1}{2}|\nabla \phi|^{2}\) :

         \begin{align*}
           \iiint_{\mathcal{V}} \Big[- \frac{p_{x}}{\rho} x  - \frac{p_{y}}{\rho} y -gz  - \frac{p_{z}}{\rho}z -3gz - \frac{3}{2}|\nabla \phi|^{2} - 3\frac{p}{\rho}  +  4 |\nabla \phi|^{2} \Big] dV -  \iint_{ \mathcal{B}} [ \phi_{t} h ]\,dx\,dy
         \end{align*}
         
         Combining the pressure and gz terms, noting \(\vec{\boldsymbol{r}} = (x,y,z)\):
         
         \begin{align*}
           \iiint_{\mathcal{V}} \Big[ - \nabla \cdot (\frac{p}{\rho}\vec{\boldsymbol{r}})  -4gz  - \frac{3}{2}|\nabla \phi|^{2}   +  4 |\nabla \phi|^{2} \Big] dV -  \iint_{ \mathcal{B}} [ \phi_{t} h ]\,dx\,dy
         \end{align*}
         
         By Divergence Theorem:
         
         \begin{align*}
           - \oiint_{\partial \mathcal{V}} \Big[ \frac{p}{\rho}\vec{\boldsymbol{r}} \Big] \cdot \vec{n}\, dS  + \iiint_{\mathcal{V}} \Big[  -4gz  - \frac{3}{2}|\nabla \phi|^{2}   +  4 |\nabla \phi|^{2} \Big] dV -  \iint_{ \mathcal{B}} [ \phi_{t} h ]\,dx\,dy
         \end{align*}
         
         Noting that pressure is 0 on the surface, and so this surface integral can reduce down to the bottom domain by decay arguments
         
         \begin{align*}
           - \iint_{\mathcal{B}} \Big[ \frac{ph}{\rho} \Big] \,dx\,dy + \iiint_{\mathcal{V}} \Big[  -4gz  - \frac{3}{2}|\nabla \phi|^{2}   +  4 |\nabla \phi|^{2} \Big] dV -  \iint_{ \mathcal{B}} [ \phi_{t} h ]\,dx\,dy
         \end{align*}
         
         Combining bottom terms, and also combining \(|\nabla \phi|^{2}\) type terms
         
         \begin{align*}
          \iiint_{\mathcal{V}} \Big[  -4gz  + \frac{5}{2}|\nabla \phi|^{2}  \Big] dV -  \iint_{ \mathcal{B}} [ \phi_{t} h +  \frac{ph}{\rho}]\,dx\,dy
         \end{align*}

         \begin{align*}
          \iiint_{\mathcal{V}} \Big[  -9gz + 5gz + \frac{5}{2}|\nabla \phi|^{2}  \Big] dV -  \iint_{ \mathcal{B}} [ \phi_{t} h +  \frac{ph}{\rho}]\,dx\,dy
         \end{align*}

         \begin{align*}
          \iiint_{\mathcal{V}} \Big[  -9gz + 5gz + \frac{5}{2}|\nabla \phi|^{2}  \Big] dV -  \iint_{ \mathcal{B}} [ \phi_{t} h +  \frac{ph}{\rho}]\,dx\,dy
         \end{align*}
         
         Integrating the \(g\eta^{2}\) terms and using Green's first Identity on the \(|\nabla \phi|^{2}\) term:

         \begin{align*}
          \iint_{\mathcal{T}} \Big[  -\frac{9}{2}g\eta^{2} + \frac{5}{2}g\eta^{2} + \frac{5}{2}q\eta_{t} \Big] \,dx\,dy - \iint_{\mathcal{B}} \Big[  -\frac{9}{2}gh^{2} + \frac{5}{2}gh^{2} \Big] \,dx\,dy -  \iint_{ \mathcal{B}} [ \phi_{t} h +  \frac{ph}{\rho}]\,dx\,dy
         \end{align*}
         
         Recognizing the Hamiltonian on the left term

         \begin{align*}
          \iint_{\mathcal{T}} \Big[  -\frac{9}{2}g\eta^{2} + 5\mathcal{H} \Big] \,dx\,dy + \iint_{ \mathcal{B}} [ 2gh^{2} - \phi_{t} h -  \frac{ph}{\rho}]\,dx\,dy
         \end{align*}
         
         But we started with (\ref{12bb}), Thus we have: 
         
         \begin{align}\label{12b}
           \frac{d}{dt}\iint_{\mathcal{T}} \Big[ (q\eta - qx\eta_{x} - qy\eta_{y})   \Big]\,dx\,dy =  \iint_{\mathcal{T}} \Big[  -\frac{9}{2}g\eta^{2} + 5\mathcal{H} \Big] \,dx\,dy + \iint_{ \mathcal{B}} [ 2gh^{2} - \phi_{t} h -  \frac{ph}{\rho}]\,dx\,dy
         \end{align}
         
         Remembering that equation (\ref{12a}) is: 
         
         \begin{align*}
           \frac{d}{dt} \iint_{\mathcal{T}} [T_{12}] \,dx\,dy =  \iint_{\mathcal{T}} \Big[ \frac{9}{2}g\eta^{2} -5\mathcal{H} + \frac{d}{dt}(q\eta - qx\eta_{x} - qy\eta_{y})   \Big]\,dx\,dy - \iint_{\mathcal{B}} \Big[ 9gt\phi\Big] \,dx\,dy
         \end{align*}
         
         Upon substituting (\ref{12b}) into (\ref{12a})
         
         \begin{align*}
           \frac{d}{dt} \iint_{\mathcal{T}} [T_{12}] \,dx\,dy =  \iint_{\mathcal{T}} \Big[ \frac{9}{2}g\eta^{2} -5\mathcal{H}   \Big]\,dx\,dy + \iint_{\mathcal{T}} \Big[  -\frac{9}{2}g\eta^{2} + 5\mathcal{H} \Big] \,dx\,dy + \iint_{ \mathcal{B}} [ 2gh^{2} - \phi_{t} h -  \frac{ph}{\rho}]\,dx\,dy - \iint_{\mathcal{B}} \Big[ 9gt\phi\Big] \,dx\,dy
         \end{align*}
         
         The top domain terms manifestly cancel, leaving: 
         
         \begin{align*}
           \frac{d}{dt} \iint_{\mathcal{T}} [T_{12}] \,dx\,dy =  \iint_{ \mathcal{B}} [ 2gh^{2} -  \frac{ph}{\rho} -  \phi_{t} h - 9gt\phi \Big] \,dx\,dy
         \end{align*}
         
         Thus, we have related the time derivative of the double iterated integral of \(T_{12}\) over the top domain to a bottom domain contribution, and thus we are done.